\documentclass[english,11pt,reqno]{amsart}

\usepackage{subfigure}
\usepackage{hyperref}
\usepackage{amssymb}
\usepackage{amsmath}
\usepackage{amsthm}
\usepackage{dsfont}
\usepackage{graphicx}
\usepackage{grffile}
\usepackage{natbib}
\usepackage{appendix}
\usepackage{grffile} 

\theoremstyle{plain}
\newtheorem{theorem}{Theorem}[section]

\newtheorem{lemma}[theorem]{Lemma}
\newtheorem{proposition}[theorem]{Proposition}

\theoremstyle{definition}
\newtheorem{definition}[theorem]{Definition}

\theoremstyle{remark}
\newtheorem{remark}{Remark}

\newtheorem{assumption}[theorem]{Assumption}

\def\di{\displaystyle}

\newcommand{\Eb}[1]{\mathbb{E}\left(#1\right)}

\newcommand{\R}{\mathbb{R}}

\newcommand{\h}{\frac{\varphi(h)}{h}}
\newcommand{\E}{\mathbb{E}}

\begin{document}

\title{A nonstandard Euler-Maruyama scheme}
\author{Fr\'ed\'eric Pierret}
\address{SYRTE UMR CNRS 8630, Observatoire de Paris, 77 avenue Denfert-Rochereau, 75014 Paris, France}


\begin{abstract}
	We construct a nonstandard finite difference numerical scheme to approximate stochastic differential equations (SDEs) using the idea of weighed step introduced by R.E. Mickens. We prove the strong convergence of our scheme under locally Lipschitz conditions of a SDE and linear growth condition. We prove the preservation of domain invariance by our scheme under a minimal condition depending on a discretization parameter and unconditionally for the expectation of the approximate solution. The results are illustrated through the geometric Brownian motion. The new scheme shows a greater behavior compared to the Euler-Maruyama scheme and balanced implicit methods which are widely used in the literature and applications.
\end{abstract}

\maketitle

\vskip 5mm
\begin{tiny}
	SYRTE UMR CNRS 8630, Observatoire de Paris, France
\end{tiny}

\tableofcontents

\section{Introduction}
The resolution of stochastic differential equations is either difficult in general or we do not have explicit solutions. Numerical schemes provide an easy way to integrate these equations but the implementation of a ``good'' integrator, in the sense of convergence order for example, is difficult (see \cite{kloeden2} and \cite{kloeden1}).\\

Despite the convergence order, a natural question on numerical schemes can be the following: \textit{Does the considered numerical scheme preserve the dynamical properties of the initial system?} \\

The usual numerical schemes, even in the deterministic case (see \cite{cresson-pierret_nsfdm} and references therein), such as Euler, Runge-Kutta and Euler-Maruyama in the stochastic case, do not preserve dynamical properties without conditions on the time-step of the numerical integration. The question is: \textit{Can we construct a stochastic numerical scheme respecting dynamical properties with a minimum of restrictions?} \\

Major work on domain invariance has been done in \cite{mil_pla_sch} and \cite{schurz} with an introduction of a class of stochastic numerical schemes, called \textit{Balanced Implicit Methods} (BIM). In \cite{schurz}, domain invariance is illustrated through multiple examples from biology, finance and marketing sciences with these methods. The domain invariance by these methods is subject to conditions on the time-step even for the numerically computed expectation. \\

The aim of this paper is to introduce the notion of nonstandard stochastic scheme based on the rules introduced in the deterministic case by R.E. Mickens (see \cite{mickens1994}, \cite{mickens2005adv}, \cite{mickens2005}). It provides a new way to create numerical schemes which preserve dynamical properties. The nonstandard rules are based on the way to construct exact numerical schemes which lead to multiple consequences to interpret the discrete derivative, integral and nonlinear terms. One of the main differences with the balanced implicit methods, is that we obtain domain invariance unconditionally for the numerically computed expectation.\\ 

The plan of the paper is as follows:\\

In Section 2, we remind classical definitions about continuous and discrete stochastic differential equations systems. In Section 3, we introduce our scheme and the assumptions made for a stochastic differential equation. In Section 4, we study the strong convergence of our scheme in the case where $f$ and $g$ are not necessarily globally Lipschitz functions. It generalizes the result obtained in \cite{higham} and \cite[Corollary 3.11 and 3.12]{hut_jen} with the Euler-Maruyama scheme in a sense of the nonstandard context. In Section 5, we study the preservation of domain invariance such as positivity which occurs in a lot of problems in scientific fields such as astronomy (see \cite{bcp}), economics, physics or more often in biology (see \cite {cps2}, \cite{cps1}). In Section 6, we illustrate numerically the scheme and its better behavior compared to the Euler-Maruyama scheme and balanced implicit methods through the geometric Brownian motion.

\section{Summary about continuous/discrete stochastic dynamical systems}

In this section, we remind classical results about continuous and discrete stochastic differential equations systems. We refer in particular to \cite{kloeden1} and \cite{oksendal} for more details and proofs. We introduce the definition of stochastic nonstandard scheme based on the rules defined by R.E. Mickens in his book \cite{mickens1994}, \cite{mickens2005adv} and \cite{mickens2005}.\\

We consider the It\^o stochastic differential equation (SDE)
\begin{equation}
\label{dyt}
dY(t)=f(Y(t))dt+g(Y(t))dW(t),\ 0\leq t\leq T,\ Y(0)\ =Y_{0},
\end{equation}
with $Y(t)\in \mathbb{R}^n$ for each $t$, $W(t)$ is a $d$-dimensional Brownian motion, $f$: $\mathbb{R}^n\rightarrow \mathbb{R}^n$ and $g$: $\mathbb{R}^n\rightarrow \mathbb{R}^{n\times d}$. \\

In many applications, the solution of the stochastic differential equation \eqref{dyt} must belong to a given domain. Such constraint is called \emph{domain invariance} and is defined as follows:

\begin{definition}
A domain $K\subset \R^m$ is said to be \emph{invariant} for the stochastic system \eqref{dyt} 
if for every initial data $Y_0\in K$ the corresponding solution $Y(t)$, $t\geq 0$, satisfies
$$P\left(\left\{Y(t)\in K,\ t\in[0,\infty[\right\}\right)=1.
$$
\end{definition}

The following theorem (see \cite[Theorem 1]{milian}) characterizes the class of functions $f$ and $g$ such that the stochastic system \eqref{dyt} preserves the domain invariance of solutions.

\begin{theorem}\label{thminv}
Let $I\subset\{1,\dots,n\}$ be a non-empty subset and $a_i,b_i\in\R$ such that $b_i>a_i$. Then, the 
set 
$$
K:=\{x\in\R^n:\ a_i\leq x_i\leq b_i,\, i\in I\}
$$ 
is invariant for the
stochastic system \eqref{dyt} if and only if
\begin{eqnarray}
	f_i (t,x) &\geq & 0 \qquad \textnormal{\emph{for}}\  x\in K\ \textnormal{\emph{such that}}\ x_i=a_i, \nonumber \\
	f_i (t,x) &\leq & 0  \qquad   \textnormal{\emph{for}}\  x\in K\ \textnormal{\emph{such that}}\ x_i=b_i,  \nonumber \\
	g_{i,j} (t,x) &= &0 \qquad  \textnormal{\emph{for}}\  x\in K\ \textnormal{\emph{such that}}\ x_i\in\{a_i,b_i\},\  j=1,\dots,d,\nonumber
\end{eqnarray}
for all $t\geq 0$ and $i\in I$. 
\end{theorem}

Let $h \in \mathbb{R}$ with $h>0$. For $k\in \mathbb{N}$, we denote by $t_k$ the discrete time defined by $t_k = kh$. 
\begin{definition}
A general one-step stochastic numerical scheme with step size $h$ and Brownian motion $W$ which computes approximations $X_k\approx Y(t_k)$ of the solution of a general system such as (\ref{dyt}) with $X_0=Y_0$ can be written in the form
\begin{equation}
	\label{sto_one_step}
	\quad X_{k+1} = \psi_{h,W}(X_k)
\end{equation}
where $\psi_{h,W}$ is a function depending on $f$ and $g$, for all $k\ge0$.
\end{definition}

\begin{remark}
In the case where $g$ is identically zero, we recover the usual definition of deterministic numerical scheme which approximate $$\frac{dY(t)}{dt}=f(Y(t)).$$ In this case $\psi_{h,W}$ is denoted only by to $\psi_{h}$.
\end{remark}

\begin{definition}
Consider a {\it continuous-time extension} $\overline{X}(t)$ of the discrete approximation \eqref{sto_one_step}. A general one-step stochastic numerical scheme $\eqref{sto_one_step}$ is said to be strongly convergent if its continuous time approximation $\overline{X}$ satisfies
\begin{equation}
	\lim_{h\rightarrow 0}\E\left[\sup_{0\leq t\leq T}|\overline{X}(t)-Y(t)|\right]=0.
\end{equation}
\end{definition}
\noindent As the continuous case, we define the domain invariance for a numerical scheme.
\begin{definition}
A domain $K$ is said to be invariant for a general one-step stochastic numerical scheme $\eqref{sto_one_step}$ if for any initial condition $X_0\in K$, $X_k$ satisfies $$P\left(\left\{X_k\in K,\ k\geq 0\right\}\right)=1.$$
\end{definition}

In the deterministic case, the rules defined by R.E. Mickens (see \cite{mickens1994}, \cite{mickens2005adv}, \cite{mickens2005}) can be states as follows:

\begin{definition}
\label{nsdef}
A general one-step deterministic numerical scheme $\eqref{sto_one_step}$ is called nonstandard finite difference scheme if at least one of the following conditions is satisfied:
\begin{itemize}
	\item $\frac{d}{dt}Y(t_k)$ is approximate as $\frac{X_{k+1}-X_k}{\varphi(h)}$ or equivalently $\int_{t_k}^{t_{k+1}}f(Y(s))ds$ is approximate as $f(X_k)\varphi(h)$ where $\varphi(h)=h+O(h^2)$ is a nonnegative function,
	\item $\psi_h (X_k)=\tilde{\psi}_h (X_{k-1},X_k,X_{k+1})$ is a nonlocal approximation of $f(Y(t_k))$.
\end{itemize}
\end{definition}

The terminology of {\it nonlocal approximation} comes from the fact that the approximation of a given function $f$ is not only given at point $X_k$ by $f(X_k)$ but can eventually depend on more points, as for example 
\begin{align*}
x^2(t_k) &\approx X_k X_{k+1}, \ X_k X_{k-1}, \ X_k \left(\frac{X_{k-1}+X_{k+1}}{2}\right), \\
x^3(t_k) &\approx X_k^2 X_{k+1}, \ X_{k-1} X_k X_{k+1}. \\
\end{align*}

\begin{remark}
In the previous definition, we restricted our attention to the easiest case, depending only on $X_{k-1}$,$ X_k$ and $X_{k+1}$. Of course, more points can be considered if necessary.
\end{remark}

\section{Nonstandard-Euler-Maruyama scheme and assumptions}
Following the first rule in Definition \ref{nsdef}, we introduce the NonStandard-Euler-Maruyama (NSEM) scheme applied to \eqref{dyt} which is given by
\begin{equation}
\label{Xk}
X_{k+1}=X_{k}+f(X_{k})\varphi(h)+g(X_{k})\Delta W_{k},
\end{equation}
where $\varphi(h)=h+c(h)$ is a nonnegative function with $c(h)=O(h^2)$ and $\Delta W_{k}=W(t_{k+1})-W(t_{k})$. A continuous-time extension $\overline{X}(t)$ of \eqref{Xk} is given by
\begin{equation}
\label{cte}
\overline{X}(t)\ :=X_{0}+\int_{0}^{t}\frac{\varphi(h)}{h}f(X(s))ds+\int_{0}^{t}g(X(s))dW(s),
\end{equation}
where $X(t)$ is defined by $X(t):=X_{k}\ \text{for}\ t\in[t_{k},\ t_{k+1}[$ . \\

Up to our knowledge, the NSEM scheme is a new numerical scheme. It can not be obtained as a specialization of a BIM scheme (see Remark \ref{remark_diff}) or a modified version of the classical EM scheme found in the literature.

\begin{remark}
\label{remark_diff}
The approach to construct balanced implicit methods (see \cite[Section 3 p.1014]{mil_pla_sch}) is based on the introduction of an explicit and implicit part in the Euler-Maruyama scheme. For a one dimensional problem, these two parts are governed by two constants such that we have
\begin{equation}
\label{bim}
X_{k+1}=X_{k}+f(X_{k})h+g(X_{k})\Delta W_{k}+(c^0 h+c^1 |\Delta W_{k}|)(X_k-X_{k+1}),
\end{equation}
where $c^0$ and $c^1$ are positive constants. For example, the BIM numerical scheme associate to the geometric Brownian motion (see Section 6) is given by (see \cite[Lemma 4.2 p.8]{schurz})
\begin{equation}
\label{bim_gbm}
X_{k+1}=\frac{X_{k}+(c^0+\lambda) X_{k}h+\sigma X_{k}\Delta W_{k}+ c^1 X_k|\Delta W_{k}|}{1+c^0h+c^1|\Delta W_{k}|},
\end{equation}
and the NSEM scheme associate to the geometric Brownian motion is given by
\begin{equation}
\label{nsem_gbm}
X_{k+1}=X_{k}+\lambda X_{k}\varphi(h)+\sigma X_{k}\Delta W_{k}.
\end{equation}
We can see that the construction is completely different. The choice of $\varphi(h)$ is given in Section 6 with a comparison of these two schemes in Section 6.2. It shows a better behavior of the NSEM scheme compared to the BIM scheme on this example.
\end{remark}

\section{Convergence of the Nonstandard-Euler-Maruyama scheme}

As the NSEM scheme can not be derived from a  known numerical scheme for stochastic differential equations (see Remark \ref{remark_diff}), we can not deduce the convergence of our scheme using classical proofs (see Remark \ref{remark_diff2} for more details). In order to study the strong convergence of the NSEM scheme, we make assumptions on the stochastic differential equations \eqref{dyt}:

\begin{assumption}[Initial condition]
\label{assum1}
We assume that the initial condition is chosen independently of the Brownian motion $W$ driving the equation \eqref{dyt}.
\end{assumption}

\begin{assumption}[Local Lipschitz condition]
\label{assum2}
For each $R>0$ there exists a constant $L_R$, depending only on $R$, such that
\begin{equation*}
	|f(x)-f(y)|\vee|g(x)-g(y)|\leq L_{R}|x-y|
\end{equation*}
for all $x,y\in \mathbb{R}^{n}$ with $|x|\vee|y|\leq R$.
\end{assumption}

\begin{assumption}[Linear Growth Condition]
\label{assum3}
There exist $K>0$ such that 
\begin{equation*}
	|f(x)|\vee|g(x)|\leq K(1+|x|)
\end{equation*}
for all $x\in\mathbb{R}^{n}$.
\end{assumption}

We follow the strategy used for the Euler-Maruyama scheme in the proof of \cite[Theorem 2.2]{higham} to prove the strong convergence of the NSEM scheme. From Assumptions \ref{assum1}, \ref{assum2} and \ref{assum3}, we deduce two results:

\begin{lemma}
\label{lemass1}
Under Assumptions \ref{assum1} and \ref{assum3}, for any $p\geq 2$ there exists a constant $C_1$ depending on $h$ and $p$ such that the exact solution and the $NSEM$ approximate solution to the Equation \eqref{dyt} have the property
\begin{equation*}
	E\left[\sup_{0\leq t\leq T}|Y(t)|^{p}\right]\vee E\left[\sup_{0\leq t\leq T}|\overline{X}(t)|^{p}\right]\leq C_1(h,p)\ .
\end{equation*}
\end{lemma}

The proof is given in Appendix A.1

\begin{lemma}
\label{lemass2}
Under Assumptions \ref{assum1}, \ref{assum2} and \ref{assum3} there exists $C_2$ depending on $h,p$ and $R$ such that
\begin{equation*}
	\E\left[\left|X(s)-\overline{X}(s)\right|^2\right] \leq C_2(h, p, R)\left(\varphi(h)^2 +nh\right)
\end{equation*}
\end{lemma}

The proof is given in Appendix A.2

\begin{remark}
In \cite{higham} the result obtained in Lemma \ref{lemass1} is an assumption for $p>2$. As noticed by the authors, it is a strong assumption but it can be recovered by assuming linear growth condition.
\end{remark}

\begin{theorem}
\label{convnsem}
Under Assumptions \ref{assum1}, \ref{assum2} and \ref{assum3}, the $NSEM$ solution of \eqref{dyt} with continuous-time extension \eqref{cte} is strongly convergent.
\end{theorem}
The proof is given in Appendix A.3

\begin{remark}
The proof of Theorem \ref{convnsem} encompasses the case where we have the global Lipschitz condition, i.e. $L_R \leq L$ for all $R$. Bounding $C_1$, $C_2$ such as they are independent of $R$ and $h$, we obtain $C_3$ which depend only on $h$ as $C_3=A\left(\varphi(h)^2 +h+\left(\frac{c(h)}{h}\right)^2\right)$ where $A$ is a constant. Then, choosing $\delta=h$ and $R=\frac{1}{h ^{1/(p-2)}}$ we obtain 
\begin{equation*}
\E\left[\sup_{0\leq t\leq T}|e(t)|^{2}\right]\leq \alpha\left(\varphi(h)^2 +h+\left(\frac{c(h)}{h}\right)^2\right)
\end{equation*}
where $\alpha$ is a constant. Moreover, in the standard case, i.e. $c(h)=0$, we recover the classical result for the Euler-Maruyama scheme,
\begin{equation*}
\E\left[\sup_{0\leq t\leq T}|\bar{X}(t)-Y(t)|^{2}\right]=O(h)
\end{equation*}
found for example in \cite{hut_jen}, \cite{kloeden1}, and \cite{mao}.
\end{remark}

\begin{remark}
\label{remark_diff2}
Our proof is not a particular case of the proof of strong convergence by local Lipschitz assumption and moment bounds for the Euler-Maruyama or the balanced implicit methods. Indeed, our scheme does not enter in the class considered in \cite[Corollary 3.11, 3.12 and 3.15]{hut_jen} and \cite[Theorem 1]{mil_pla_sch} by construction. However, different hypothesis to prove the strong convergence of the Euler-Maruyama scheme can be used to weaken the linear growth condition. For a complete review and details of these possibilities for the Euler-Maruyama scheme, see \cite[Section 3.4]{hut_jen}.
\end{remark}

\section{Domain invariance}
In order to study the invariance describe in Theorem \ref{thminv} for the numerical approximation of $\eqref{dyt}$ obtained by the NSEM scheme, we restrict our intention to the domain $K^+:=\{x\in\R^n : \ a_i\leq x_i,\ i\in I \}$ and the domain $K^-:=\{x\in\R^n:\ x_i \leq b_i, \ i\in I\}$ where $I\subset\{1,\dots,n\}$ is a non-empty subset and $a_i,b_i\in\R$ such that $b_i>a_i$. As $K=K^+ \cap K^-$ the condition of invariance for the numerical scheme of $K$ would be both conditions of invariance for $K^+$ and $K^-$. The methods being the same to prove the invariance of $K^+$ and $K^-$, we write only the proof for $K^+$. \\

As there exists a multitude of choices for the function $\varphi$, we select one of the most used in the literature of nonstandard schemes, which is the following:
\begin{equation}
	\label{nsstep}
	\varphi(h)=\frac{1-\phi(\alpha h)}{\alpha}
\end{equation}
where $\alpha$ is a positive constant and $\phi$ is a function which takes values in $]0,1[$. \\

We consider the case where $g$ is identically zero, that is to say the NSEM scheme reduces to the nonstandard Euler scheme (see \cite{mickens1994}). The domain invariance in that case is well-known in the literature but we remind it in order to show the steps to prove the nonstandard case. We have:

\begin{proposition}
	\label{dominv1}
	Let $K$ be an invariant domain for $Y(t)$ and
	\begin{equation*}
		D=\sup_{x\in K} \left\|\int_0^1 f'(t_k,a+s(x-a))ds\right\|_\infty.
	\end{equation*}
	For all $h>0$, $K$ is invariant for the numerical approximation of $Y(t)$ obtained by the nonstandard Euler scheme. The positive constant $\alpha$ in \eqref{nsstep} is chosen to be $D$.
\end{proposition}
The proof is given in Appendix \ref{dem_dominv1} \\

\begin{remark}
	We arbitrary choose the uniform matrix norm $\|\cdot \|_\infty$ for $f'$.
\end{remark}

Now, we study the invariance of $K$ for the numerical approximation of $\eqref{dyt}$ obtained by the NSEM scheme in the general case. There is a huge difference between the classical case ($g$ is identically zero) and the stochastic case ($g$ is not identically zero) due to the fact that $\Delta W_k$ is almost surely unbounded (see \cite{oksendal}) for all $k\ge 0$. However, we can estimate the probability such that the NSEM scheme will respect the invariance of $K$. \\

\begin{proposition}
	\label{propcondinv}
	Let $K$ be an invariant domain for $Y(t)$ and 
	\begin{equation*}
		D=\sup_{x\in K} \left\|\int_0^1f'(t_k,a+s(x-a))ds\right\|_\infty, \ S=\sup_{x\in K} \left\|\int_0^1 g'(a+s(x-a))ds\right\|_\infty.
	\end{equation*}
	If $h>0$ is such that
	\begin{equation*}
		\sup_{1\leq p \leq d}|\Delta W_k|_p\leq \frac{\phi(Dh)}{Sd}
	\end{equation*}
	for all $k\ge 0$  then, the domain $K$ is invariant for the numerical approximation of $Y(t)$ obtained by the NSEM scheme.
\end{proposition}
The proof is given in Appendix \ref{dem_dominv1_nsem}. \\

\begin{remark}
	Again, we arbitrary choose the uniform norm $\|\cdot \|_\infty$ for $f'$ and $g'$. Let us remark that $g'$ is ``cubic matrix'' or more precisely, a tensor of order 3.
\end{remark}
We quantify the probability that the condition
\begin{equation*}
	\sup_{1\leq p \leq d}|\Delta W_k|_p\leq \frac{\phi(Dh)}{Sd}
\end{equation*}
occurs for all $k\geq0$:

\begin{proposition}
	\label{dominv2}
	For all $0<\epsilon<\frac{1}{2}$, there exists $h_0$ such that for all $h<h_0(\epsilon)$ then
	\begin{equation*}
		\mathbb{P}\left(\sup_{1\leq p \leq d}|\Delta W_k|_p\leq \frac{\phi(Dh)}{Sd} \right)>1-\epsilon
	\end{equation*}
	for all $k\ge 0$.
\end{proposition}
The proof is given in Appendix \ref{dem_dominv2}.\\

As a consequence, if the time-step is sufficiently small then, we obtain domain invariance of $K$ for the numerical approximation of $Y(t)$, obtained by the NSEM scheme. It is obtained with a probability as close as we want to one, in choosing a $\epsilon$ which represent a tolerance parameter.

\begin{proposition}
	\label{dominv3}
	Let $K$ be an invariant domain for $Y(t)$ and
	\begin{equation*}
		D=\sup_{x\in K} \left\|\int_0^1f'(t_k,a+s(x-a))ds\right\|_\infty.
	\end{equation*}
	For all $h>0$, $K$ is invariant for the numerical approximation of $\Eb{Y(t)}$ obtained by the NSEM scheme. The positive constant $\alpha$ in \eqref{nsstep} is chosen to be $D$.
\end{proposition}

The proof is given in Appendix \ref{dem_dominv3}

\section{Numerical illustrations}

We consider the geometric Brownian motion which is driven by the following stochastic differential equations:
\begin{equation}
	\label{dytex}
	dY(t)=\mu Y(t) dt+\sigma Y(t)dW(t),
\end{equation}
where $\mu$ and $\sigma$ are real constants. Without loss of generality, we assume that $\sigma$ is positive. The solution of \ref{dytex} is given by (see \cite{oksendal})
\begin{equation}
	Y(t)=Y_0\exp\left(\left(\mu - \frac{\sigma^2}{2}\right)t+\sigma W(t)\right).
\end{equation}
and its expectation is given by
\begin{equation}
	\Eb{Y(t)}=Y_0\exp\left(\mu t\right).
\end{equation}

\noindent{\textbf{Example 1}} \\

We consider the geometric Brownian motion with a negative value on its deterministic part (also called the stochastic decay equation). Let $\mu=-\lambda$ where $\lambda$ is strictly positive. In this example, the constants $D$ and $S$ of the previous section, correspond respectively to $\lambda$ and $\sigma$. The classical decay equation (the geometric Brownian motion with $\sigma =0$) has been studied by R.E. Mickens with the nonstandard Euler scheme using $\phi(\lambda h)=\exp(-\lambda h)$. We keep this choice for our example. \\

It is clear that the solution of the geometric Brownian motion equation belongs to $K^+$. Using the Proposition \ref{propcondinv}, we obtain:
\begin{lemma}
	\label{lemmeh0em}
	Let $0<\epsilon<\frac{1}{2}$. The Euler-Maruyama scheme preserves the positivity of the geometric Brownian motion equation if
	\begin{equation*}
		\di h<h_0(\epsilon)=\frac{1}{\lambda }  + \frac{ \alpha(\epsilon) ^2  \sigma ^2}{\lambda ^2} - \frac{\alpha \sigma  \sqrt{ \alpha(\epsilon) ^2  \sigma ^2+2 \lambda }}{\lambda ^2}
	\end{equation*}
	where $\alpha(\epsilon)=\mathsf{erf}^{-1}(1-\epsilon)$.
\end{lemma}
The proof is given in Appendix C.1. \\

In a same way, we have for the NSEM:
\begin{lemma}
	\label{lemmeh0ns}
	Let $0<\epsilon<\frac{1}{2}$. The Nonstandard-Euler-Maruyama scheme preserves the positivity of the geometric Brownian motion if
	\begin{equation}
		h<h_0(\epsilon)\equiv\frac{1}{2\lambda}\mathcal{W}\left(\frac{\lambda}{ \sigma^2  \alpha(\epsilon)^2}\right)
	\end{equation}
	where $\mathcal{W}$ is the product logarithm function and $\alpha(\epsilon)=\mathsf{erf}^{-1}(1-\epsilon)$.
\end{lemma}
The proof is given in Appendix C.2. \\

Choosing $c^0$ and $c^1$ respectively as $\lambda$ and $\sigma$, we obtain from \cite[Lemma 4.2 p.8]{schurz} that the BIM scheme \eqref{bim_gbm} preserves the positivity for all $h>0$ . \\ 

\begin{figure}[ht!]
		\begin{tabular}{ll}
			\subfigure[$h=T/256$]{\includegraphics[height=0.5\textwidth]{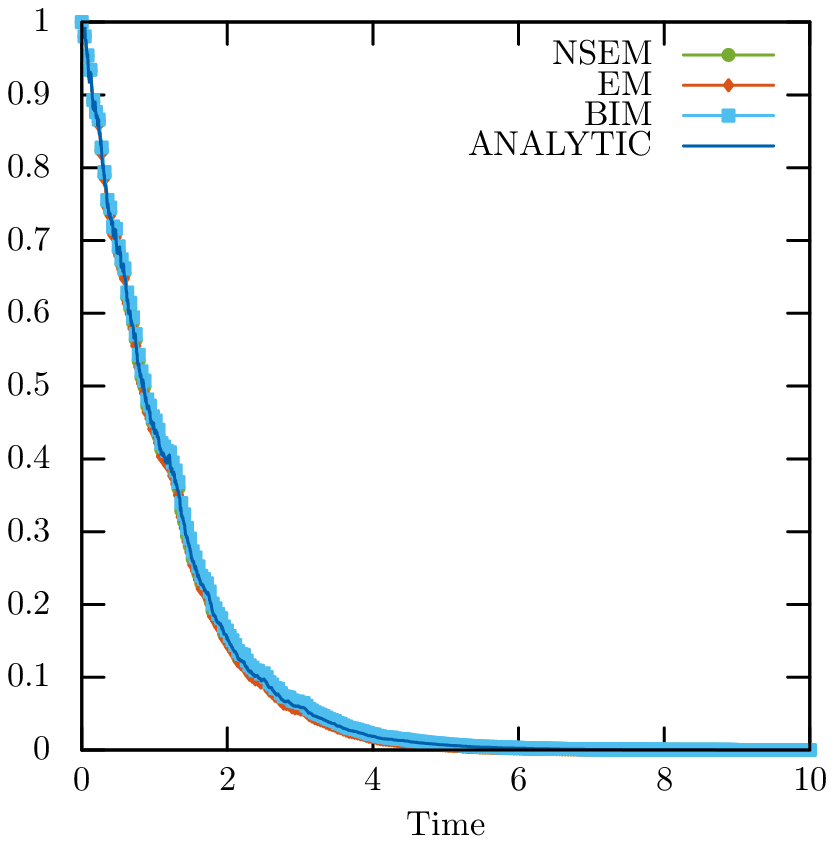}\label{a1}}&
			\subfigure[$h=T/64$]{\includegraphics[height=0.5\textwidth]{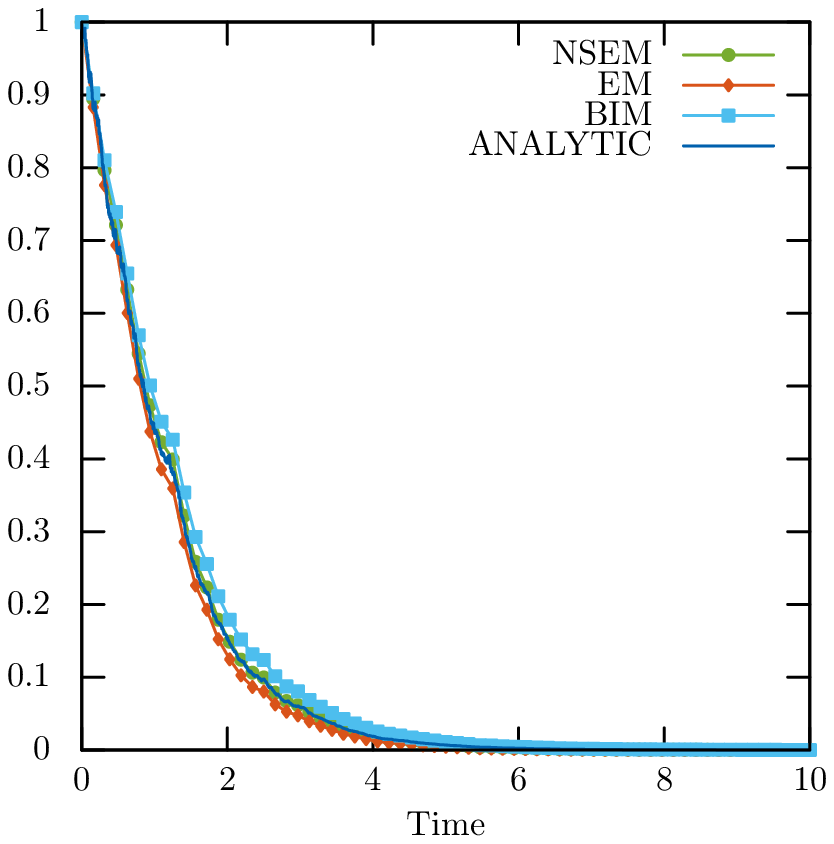}\label{b1}}\\
			\subfigure[$h=T/32$]{\includegraphics[height=0.5\textwidth]{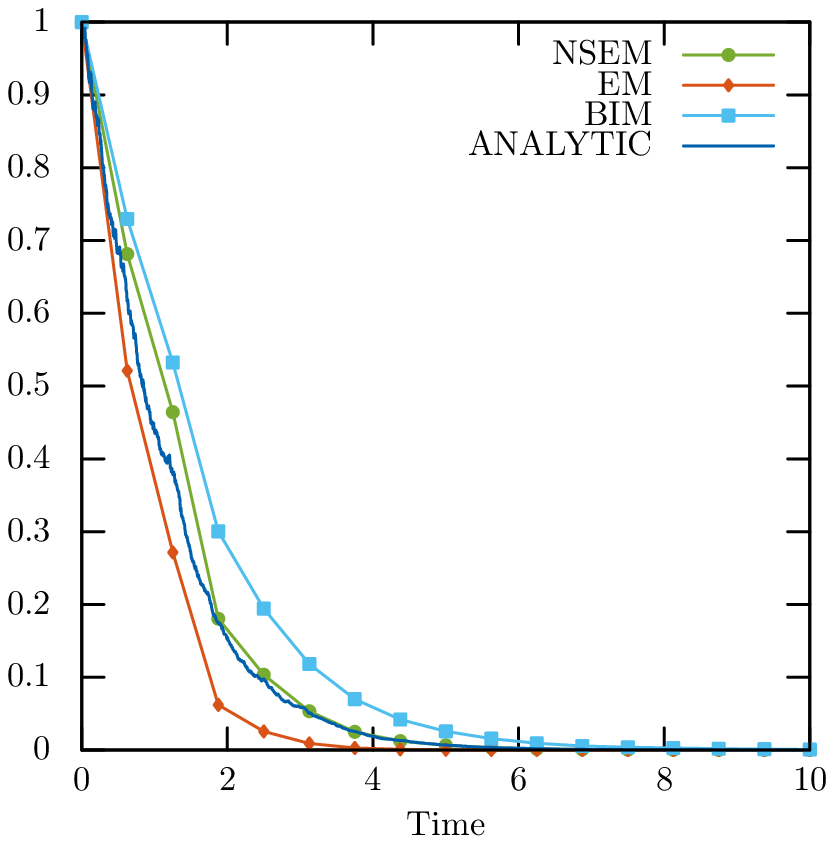}\label{c1}}&
			\subfigure[$h=T/4$]{\includegraphics[height=0.5\textwidth]{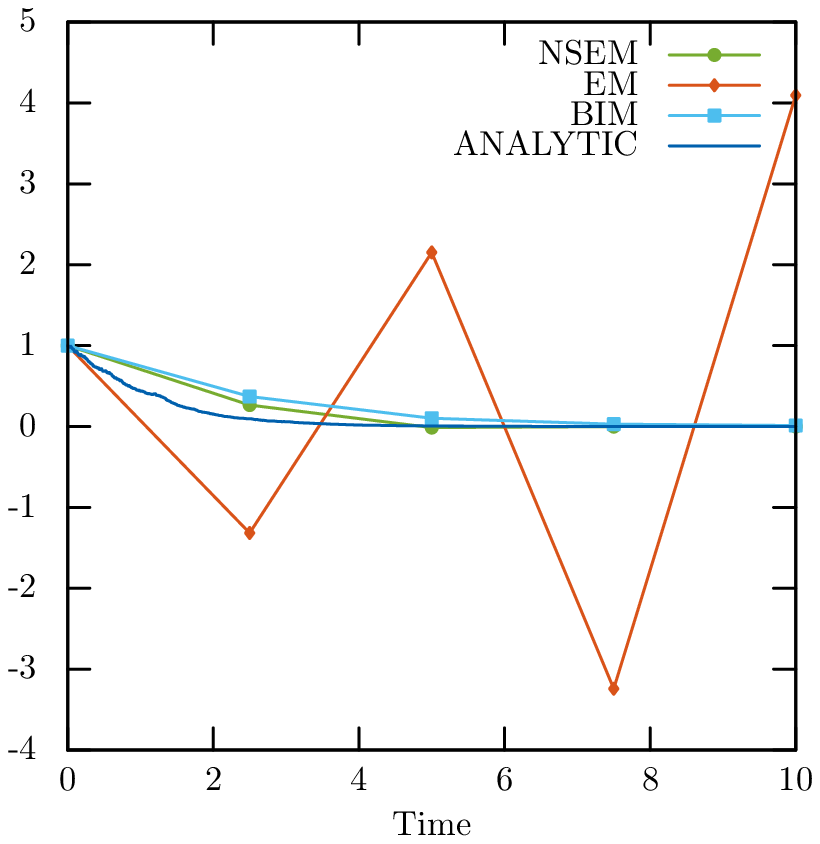}\label{d1}}\\
		\end{tabular}
	\caption{Numerical simulations of \eqref{dytex} with $Y_0=1$, $\mu=-\lambda=-1$ and $\sigma=0.1$.}
	\label{figsim1}
\end{figure}

\begin{figure}[ht!]
	\begin{tabular}{rr}
		\subfigure[$h=T/256$]{\includegraphics[height=0.5\textwidth]{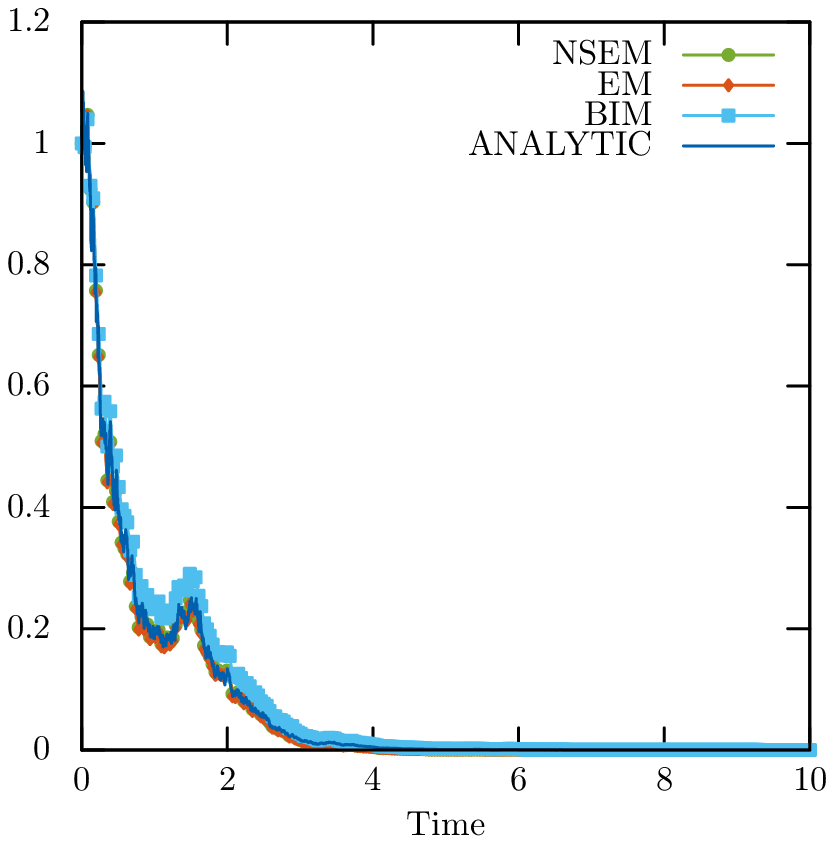}\label{a11}}&
		\subfigure[$h=T/64$]{\includegraphics[height=0.5\textwidth]{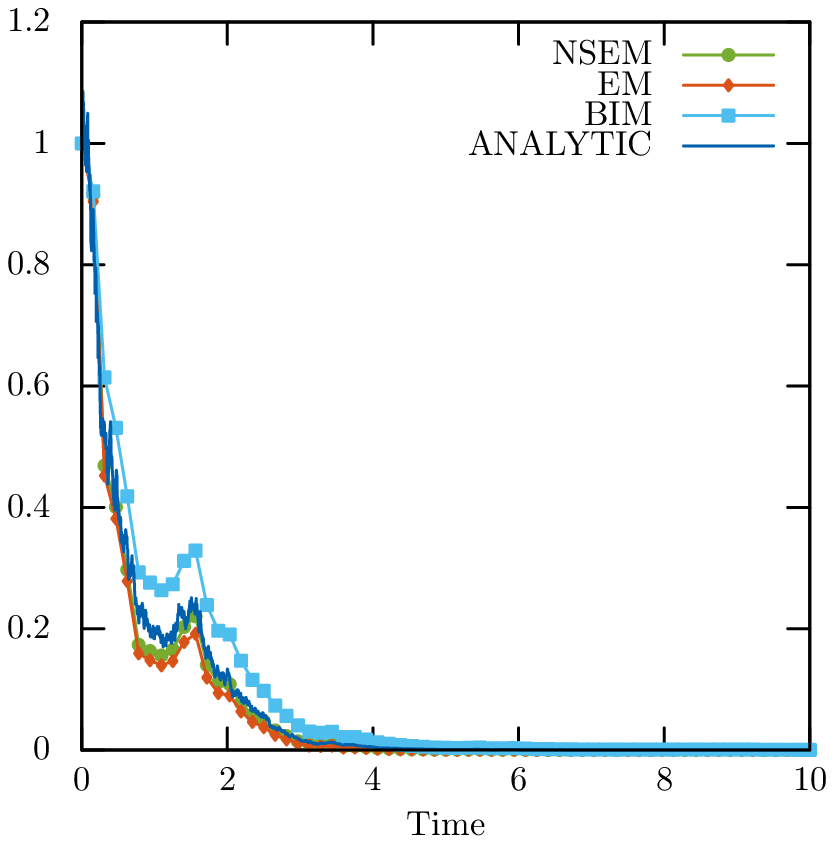}\label{b11}}\\
		\subfigure[$h=T/32$]{\includegraphics[height=0.5\textwidth]{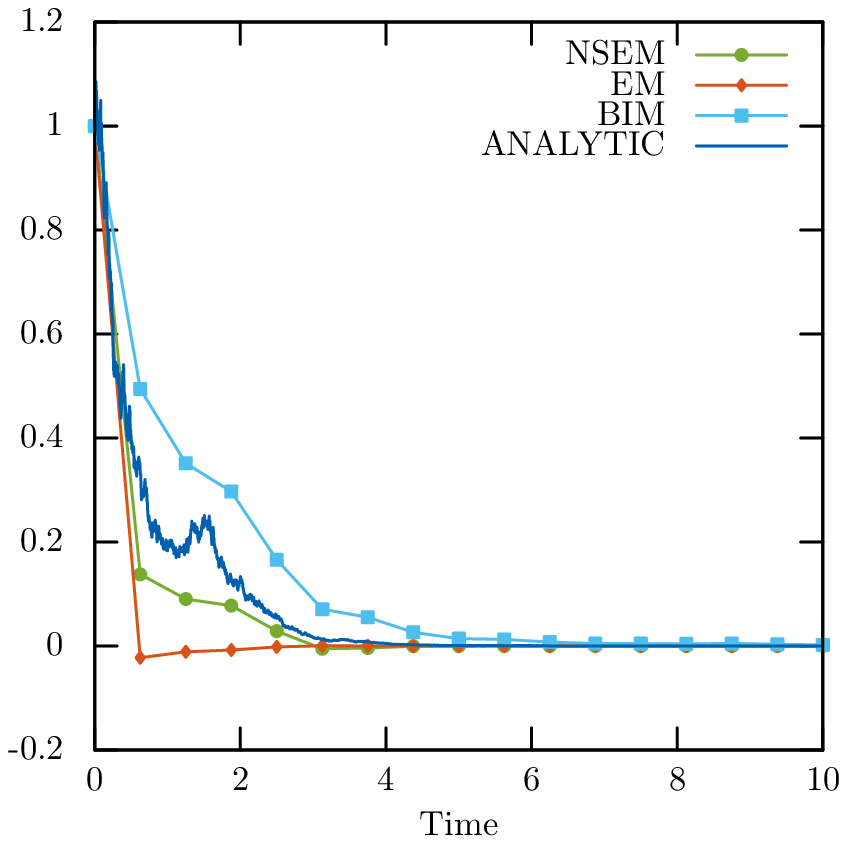}\label{c11}}&
		\subfigure[$h=T/4$]{\includegraphics[height=0.5\textwidth]{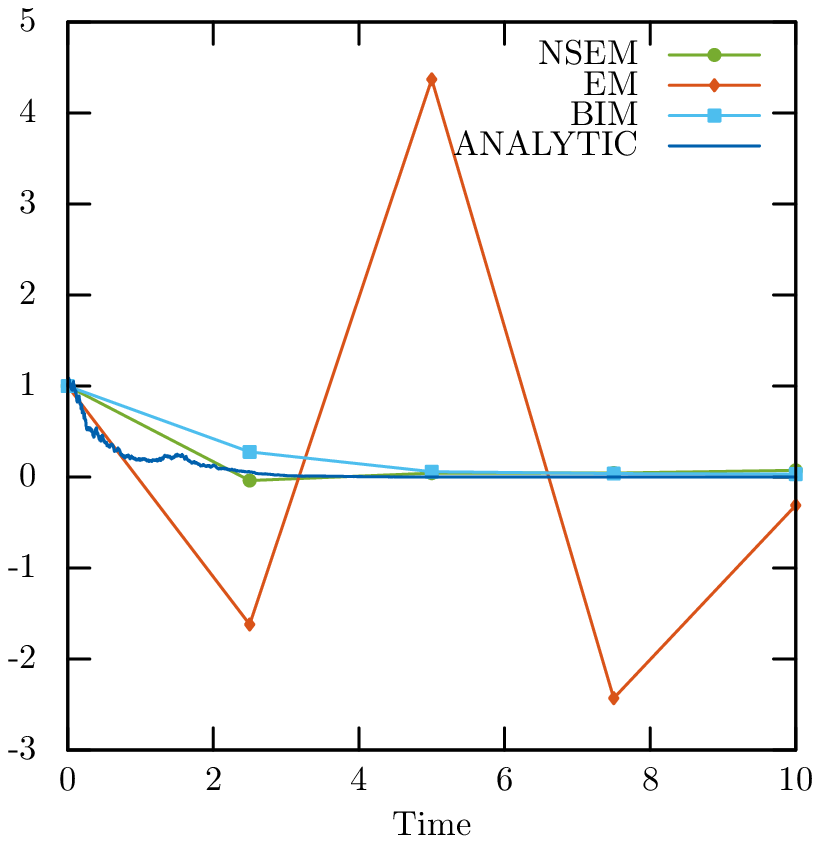}\label{d11}}\\
	\end{tabular}
	\caption{Numerical simulations of \eqref{dytex} with $Y_0=1$, $\mu=-\lambda=-1$ and $\sigma=0.5$.}
	\label{figsim11}
\end{figure}

In Figures \ref{figsim1} and \ref{figsim11}, we display simulations performed with multiple time-step over a period $T=10$ in order to show the behavior of our scheme compared to the EM scheme and the BIM scheme. The constants are chosen such that $\lambda=1$, $\sigma=0.1$ in Figure \ref{figsim1} and $\sigma=0.5$ in Figure \ref{figsim11}. \\

\begin{figure}[h!]
	\centering
	\includegraphics[width=0.5\textwidth]{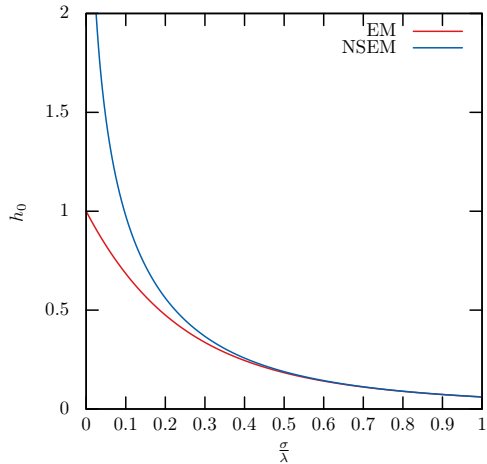}
	\caption{Minimal step for positivity invariance with respect to $\frac{\sigma}{\lambda}$ with $\epsilon=0.01$.}
	\label{figratio}
\end{figure}

Using Lemmas \ref{lemmeh0em} and \ref{lemmeh0ns}, in Figure \ref{figratio} we display the minimal step $h_0$ with respect to the ratio $\frac{\sigma}{\lambda}$ for a tolerance parameter $\epsilon=0.01$ for our scheme and EM. As we can see, the minimal step for our scheme to preserve the positivity becomes very high as long as the noise becomes smaller compared to the deterministic part. For $\lambda=1$ and $\sigma=0.1$, the minimal step for invariance by EM is $h_0^{EM}\approx 0.77$ and for NSEM, it is $h_0^{NSEM}\approx 1.25$. For $\lambda=1$ and $\sigma=0.5$, the minimal step for invariance by EM is $h_0^{EM}\approx 0.30$ and for NSEM, it is $h_0^{NSEM}\approx 0.32$. \\

In Figures \ref{figsim1} and \ref{figsim11}, we can see that the NSEM scheme globally behaves better than the EM scheme and the BIM scheme. First, in Figures \ref{a1}, \ref{a11}, \ref{b1} and \ref{b11}, we can see that the three schemes behave similarly. Knowing the minimal steps for positivity invariance, in Figures \ref{b1} and \ref{b11}, we can see that contrary to our scheme, the EM scheme does not respect the positivity as expected by Lemma \ref{lemmeh0ns}. In Figures \ref{d1} and \ref{d11}, the EM scheme diverges and the NSEM scheme still behaves correctly even if the time-step is greater than his minimal step for positivity invariance. The results in Figure \ref{d1} hold for a lot of realizations of Brownian motion because our boundaries, in the proof of minimal step for invariance in the case of the NSEM, are not optimal but give globally a very good information.\\

\begin{figure}[ht!]
	\begin{tabular}{rr}
		\subfigure[$h=0.1$]{\includegraphics[height=0.5\textwidth]{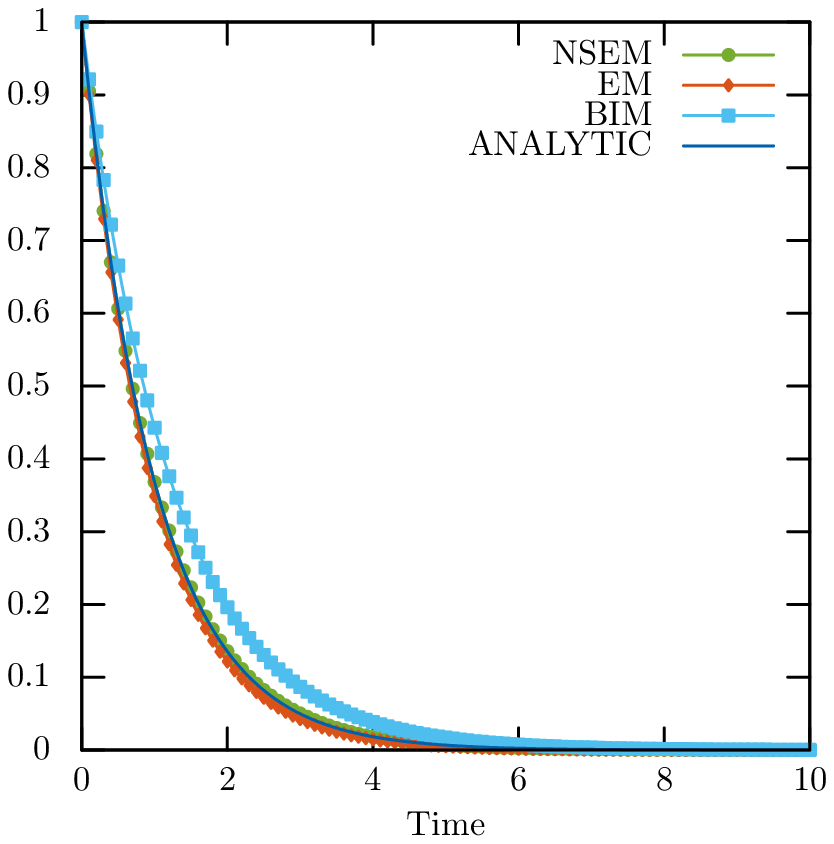}\label{a2}}&
		\subfigure[$h=1$]{\includegraphics[height=0.5\textwidth]{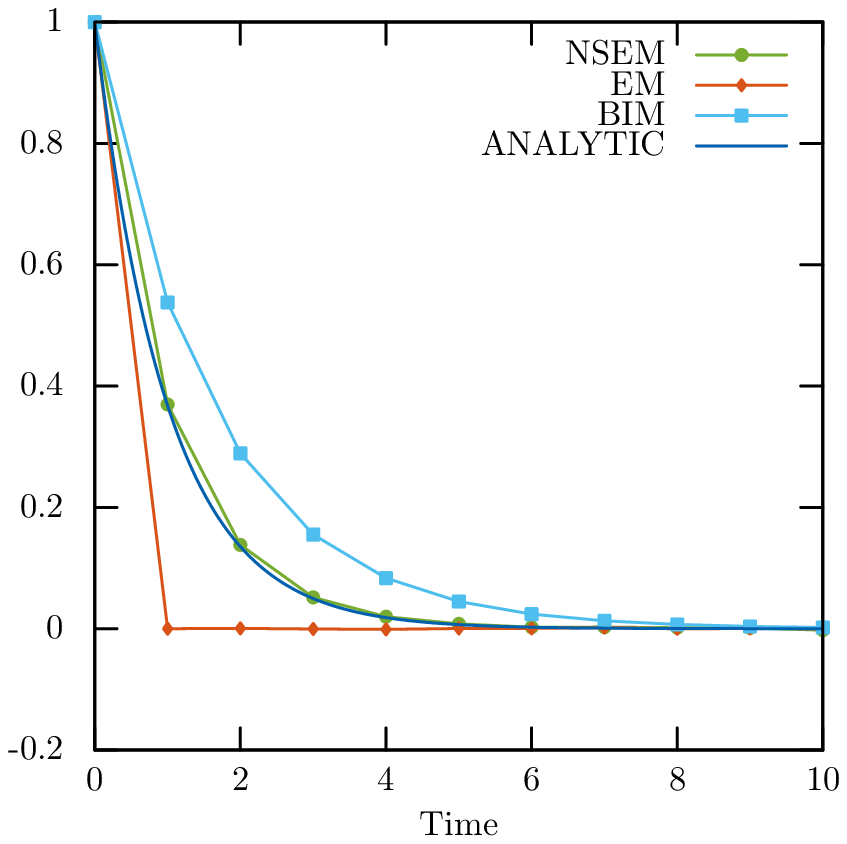}\label{b2}}\\
		\multicolumn{2}{c}{\subfigure[$h=2$]{\includegraphics[height=0.5\textwidth]{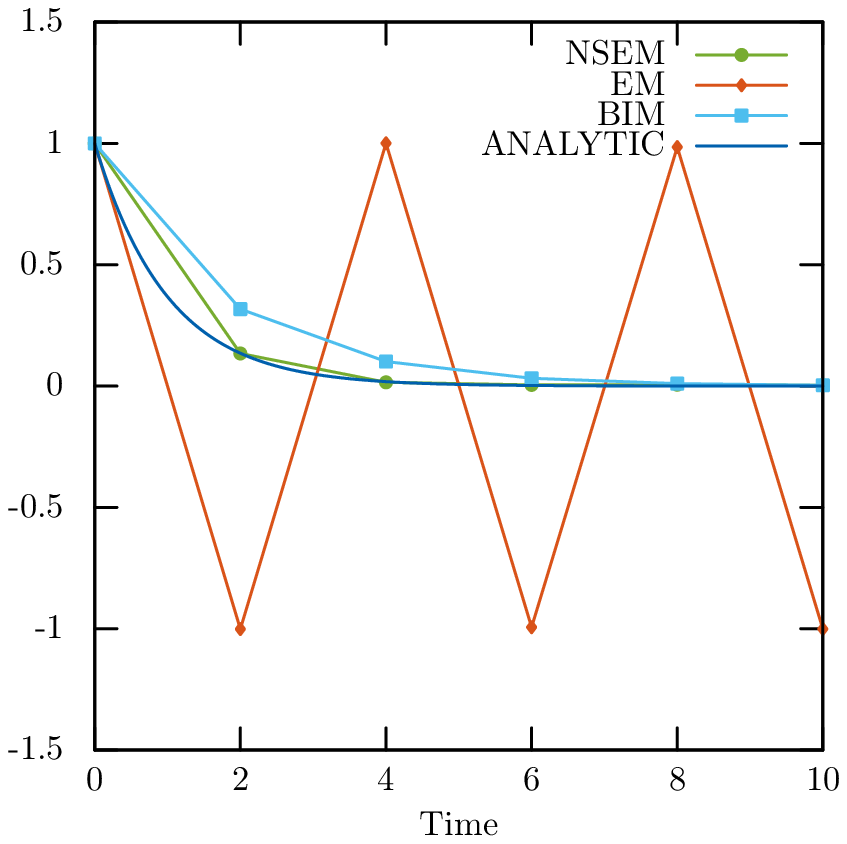}\label{c2}}}

	\end{tabular}
	\caption{Numerical simulations of the expectation of the solution of \eqref{dytex} with $Y_0=1$, $\mu=-\lambda=-1$ and $\sigma=1$.}
	\label{figsim2}
\end{figure}

In order to illustrate the advantage of our scheme when taking the expectation of the numerical solution, in Figure \ref{figsim2}, we display simulations performed with a noisy system, i.e. $\lambda=1$ and $\sigma=1$. The numerical expectation is computed using a Monte-Carlo method with $10^6$ realizations. \\

In Figure \ref{a2}, we can see that the three schemes behave similarly for $h=0.1$ with a good approximation to the expectation of the analytic solution. In Figure \ref{b2}, for $h=1$, we can see that, contrary to our scheme, the EM scheme and the BIM scheme do not approximate properly the solution from time $T=0$ to $T=5$. In Figure \ref{c2}, for $h=2$ the EM scheme oscillates around the solution whereas our scheme still behaves very well and is very closed to the analytic solution. This is exactly the case where the nonstandard character greatly improves the results contrary to the usual stochastic schemes. \\

\noindent{\textbf{Example 2}} \\

We consider the geometric Brownian motion with a positive value on its deterministic part. Let $\mu=\lambda$ where $\lambda$ is strictly positive. In this example, the constants $D$ and $S$ of the previous section, correspond respectively to $\lambda$ and $\sigma$. We still use $\phi(\lambda h)=\exp(-\lambda h)$. It is clear that the solution of the geometric Brownian equation belongs to $K^+$ and the three schemes satisfy the positivity invariance. \\

\begin{figure}[ht!]
	\begin{tabular}{rr}
		\subfigure[$h=T/256$]{\includegraphics[height=0.5\textwidth]{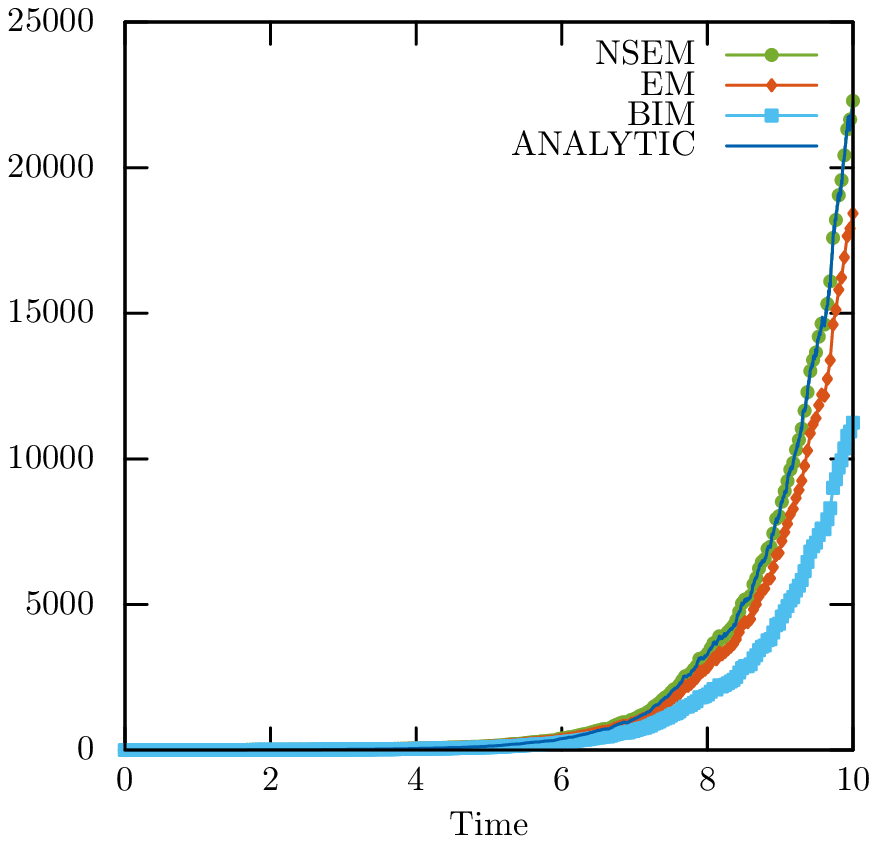}\label{a3}}&
		\subfigure[$h=T/64$]{\includegraphics[height=0.5\textwidth]{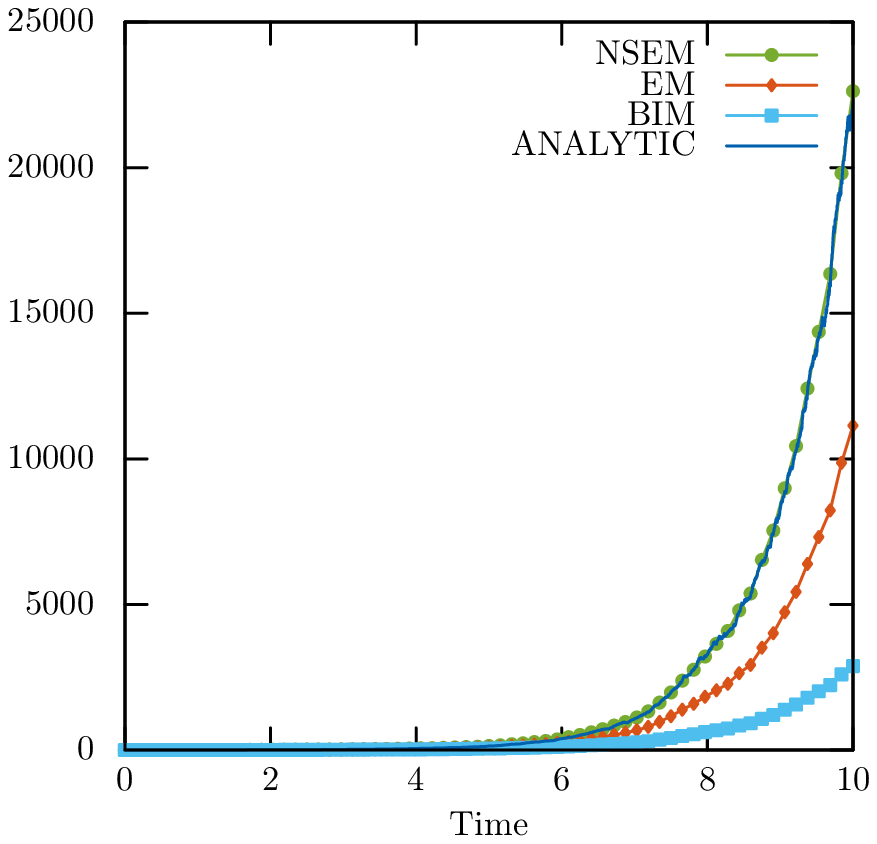}\label{b3}}\\
		\subfigure[$h=T/32$]{\includegraphics[height=0.5\textwidth]{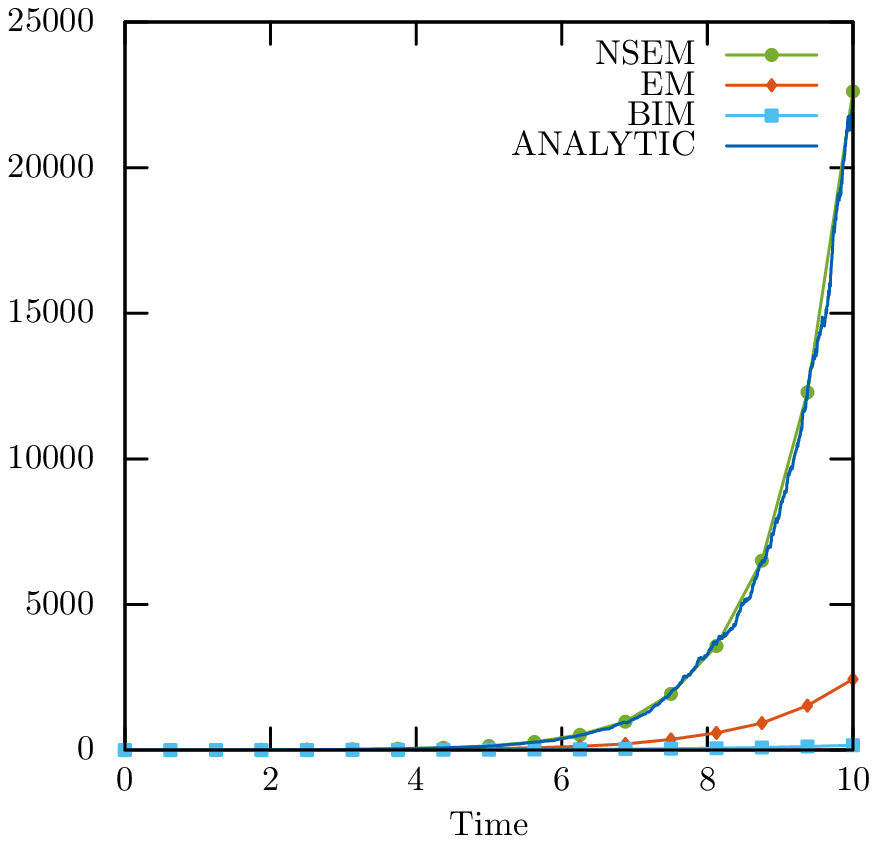}\label{c3}}&
		\subfigure[$h=T/4$]{\includegraphics[height=0.5\textwidth]{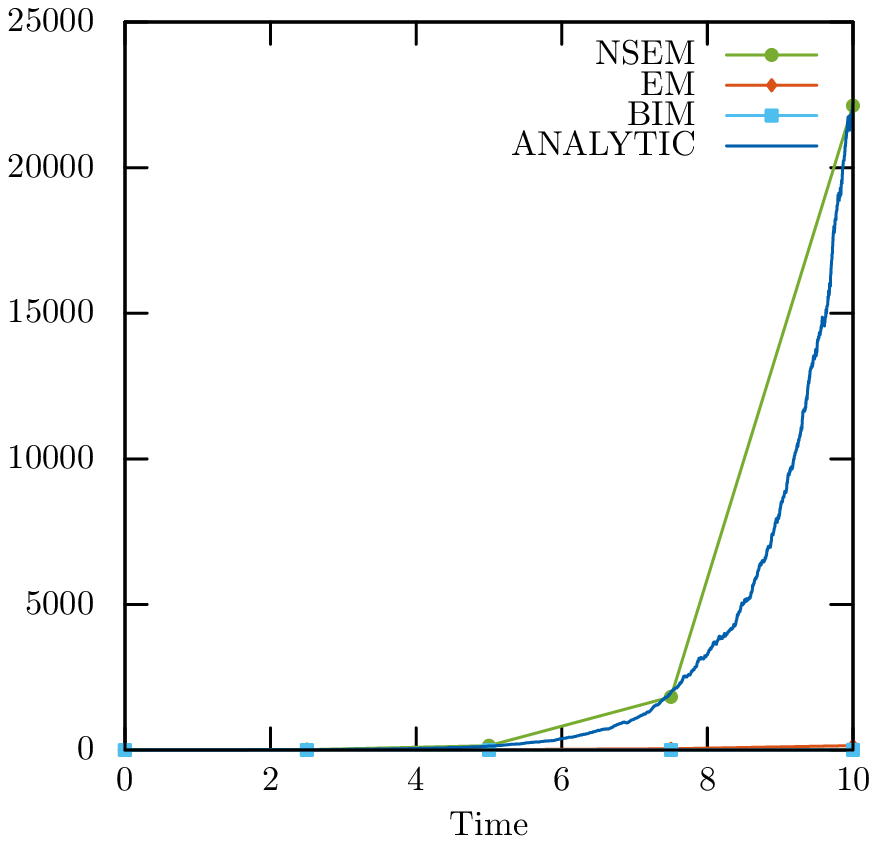}\label{d3}}\\
	\end{tabular}
	\caption{Numerical simulations of \eqref{dytex} with $Y_0=1$, $\mu=\lambda=1$ and $\sigma=0.1$.}
	\label{figsim3}
\end{figure}

\begin{figure}[ht!]
	\begin{tabular}{rr}
		\subfigure[$h=T/256$]{\includegraphics[height=0.5\textwidth]{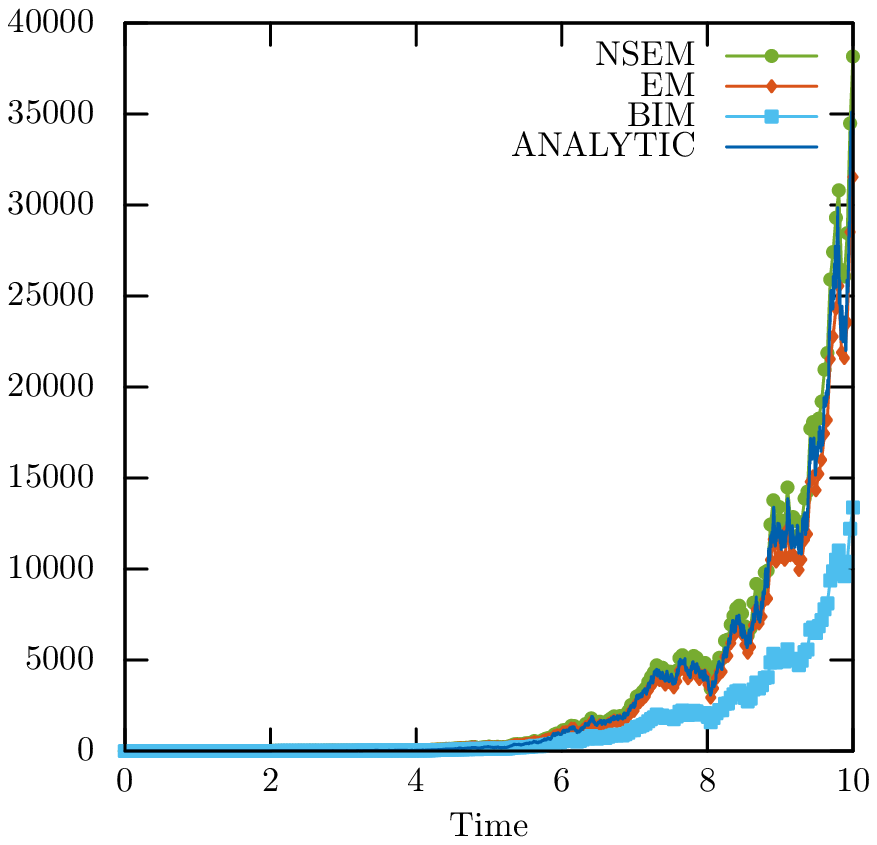}\label{a33}}&
		\subfigure[$h=T/64$]{\includegraphics[height=0.5\textwidth]{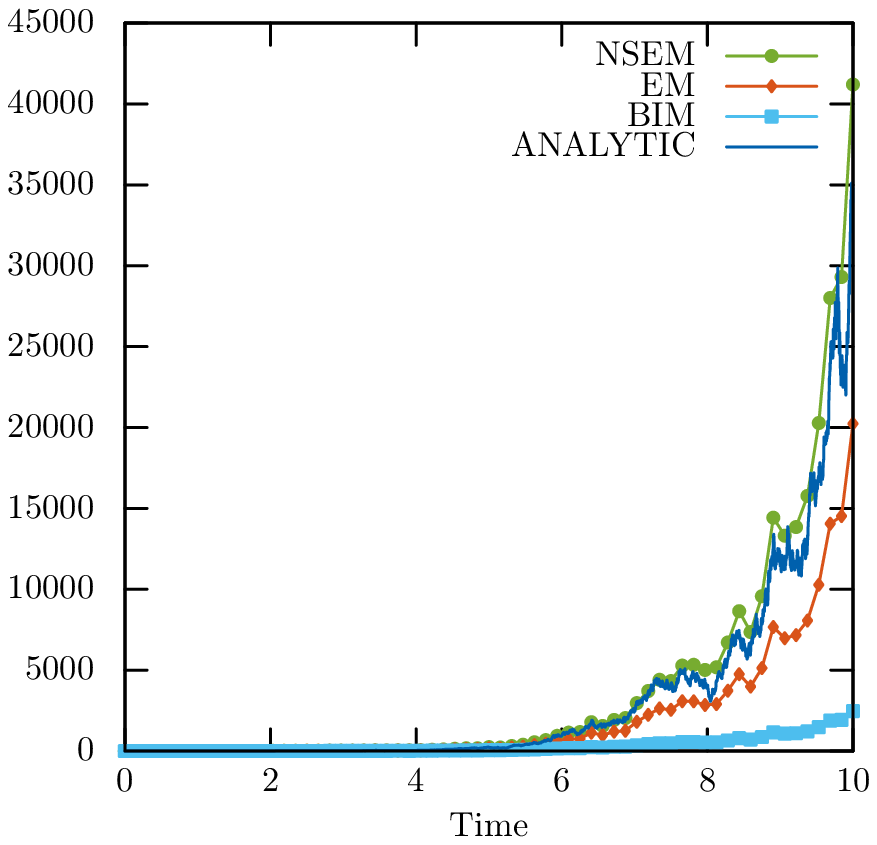}\label{b33}}\\
		\subfigure[$h=T/32$]{\includegraphics[height=0.5\textwidth]{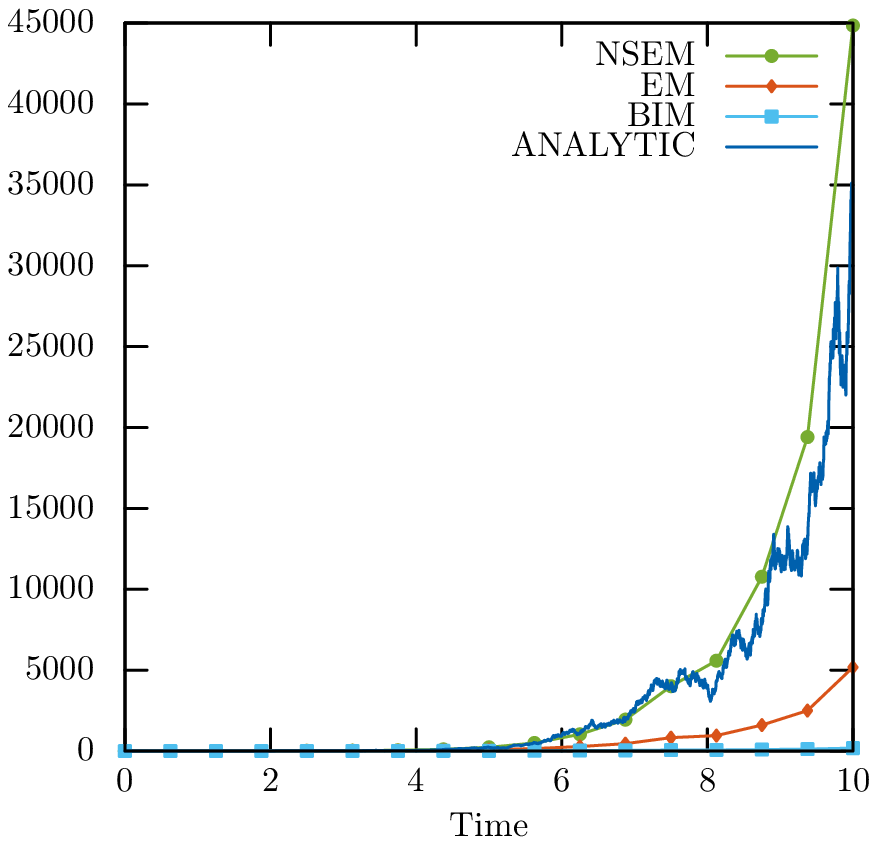}\label{c33}}&
		\subfigure[$h=T/4$]{\includegraphics[height=0.5\textwidth]{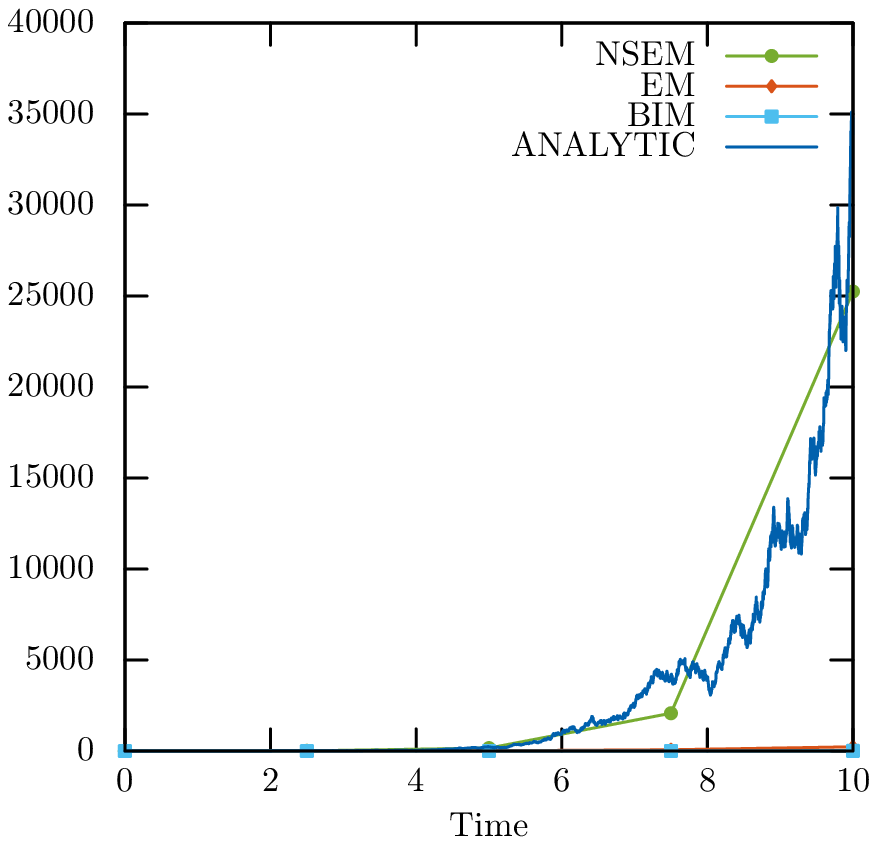}\label{d33}}\\
	\end{tabular}
	\caption{Numerical simulations of \eqref{dytex} with $Y_0=1$, $\mu=\lambda=1$ and $\sigma=0.5$.}
	\label{figsim33}
\end{figure}

As in Example 1, in Figures \ref{figsim3} and \ref{figsim33} we display simulations performed with multiple time-step over a period $T=10$ in order to show the behavior of our scheme compared to the EM scheme and the BIM scheme. The constants are chosen such that $\lambda=1$, $\sigma=0.1$ in Figure \ref{figsim3} and $\sigma=0.5$ in Figure \ref{figsim33}. \\

In Figures \ref{figsim3} and \ref{figsim33}, we can see for different time-step and values of $\sigma$ that the NSEM scheme and the BIM scheme behave similarly during the first half period integration while in the second part the BIM scheme fails to approximate the solution. \\

\begin{figure}[h!]
	\begin{tabular}{rr}
		\subfigure[$h=0.1$]{\includegraphics[height=0.5\textwidth]{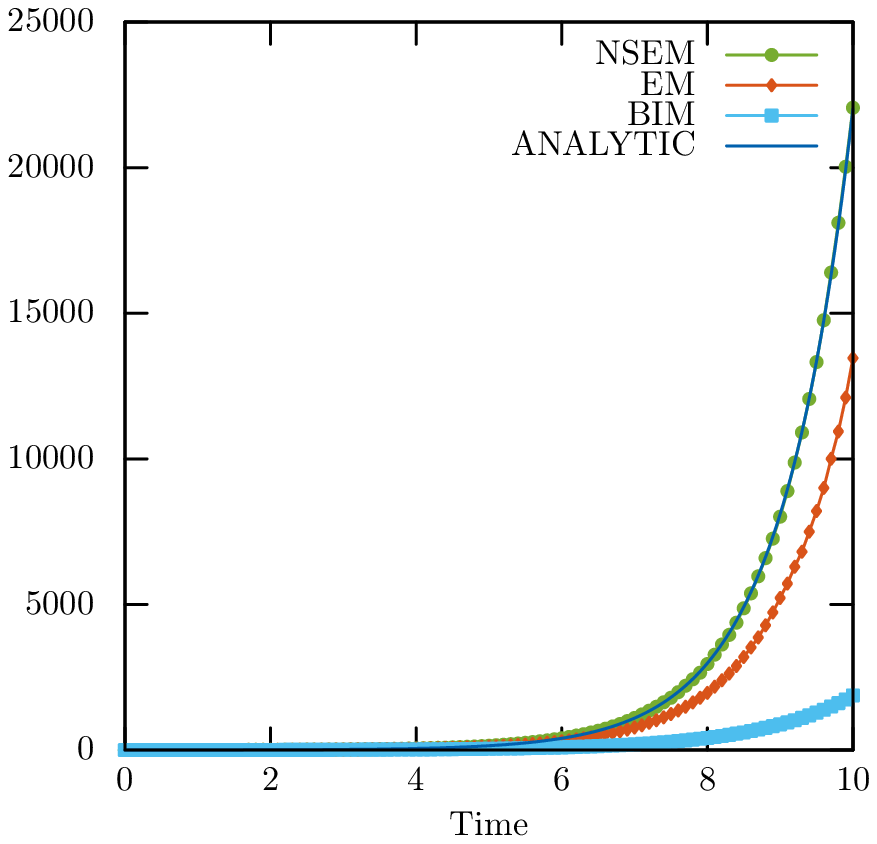}\label{b4}}&
		\subfigure[$h=1$]{\includegraphics[height=0.5\textwidth]{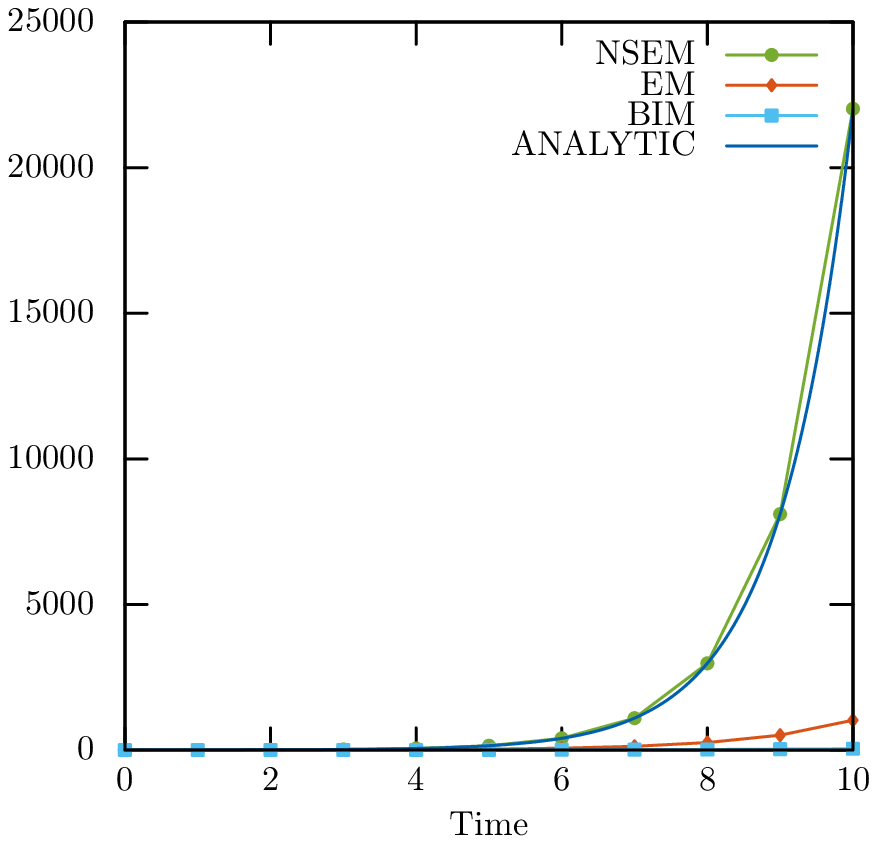}\label{c4}}\\
		\multicolumn{2}{c}{\subfigure[$h=2$]{\includegraphics[height=0.5\textwidth]{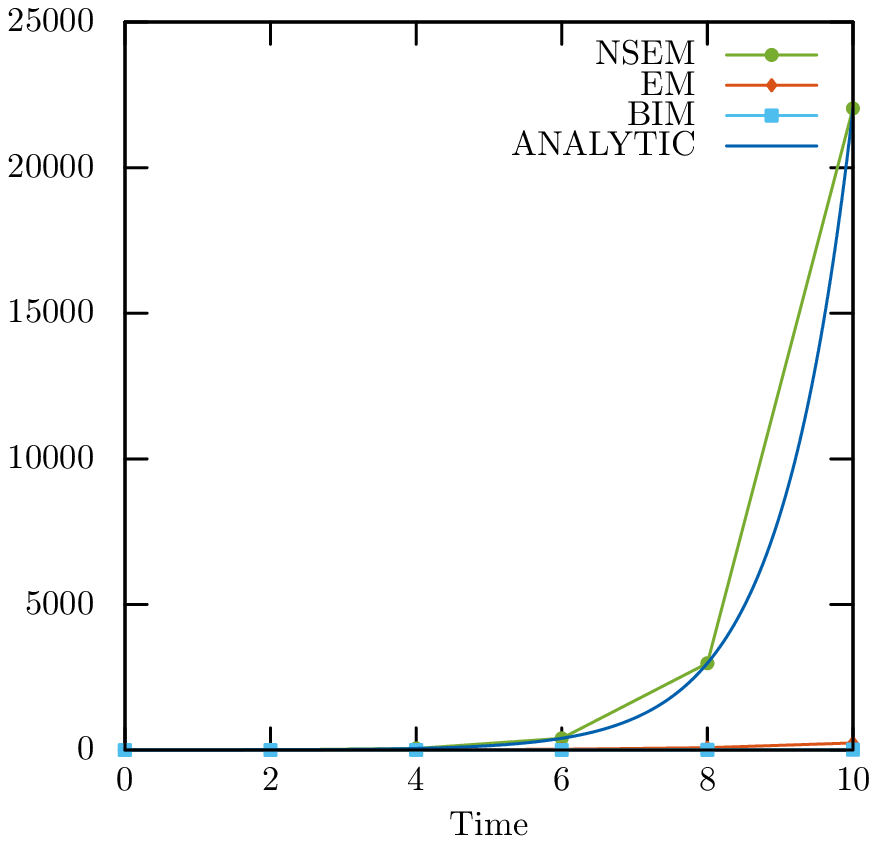}\label{d4}}}
		
	\end{tabular}
	\caption{Numerical simulations of the expectation of \eqref{dytex} with $Y_0=1$, $\mu=-\lambda=-1$ and $\sigma=1$.}
	\label{figsim4}
\end{figure}

As in Example 1, in order to illustrate again the advantage of our scheme when taking the expectation of the numerical solution, in Figure \ref{figsim4} we display simulations performed with a noisy system, i.e. $\lambda=1$ and $\sigma=1$. The numerical expectation is computed using a Monte-Carlo method with $10^6$ realizations. \\

We observe the same behavior as previous in Figure \ref{figsim2}. In all the cases, the NSEM scheme almost approximates exactly the solution due to the high number of realizations. For $h=0.1$, $h=1$ and $h=2$, the EM scheme and the BIM scheme fail completely to approximate the solution. For $h=0.1$ the EM scheme behaves greater that the BIM scheme which seem normal. Indeed, the use of the BIM scheme in such a case seems inappropriate due to the fact that by construction, the positivity is trivially satisfied. But in a more general case, the EM scheme could be more efficient and easily implementable than the BIM scheme. It gives another one advantage to our scheme which provide an efficient and easy way to tweak the EM scheme and to improve the results obtained.

\section{Conclusion and prospects}
A nonstandard Euler-Maruyama scheme has been introduced and its convergence has been proved with local Lipschitz condition and linear growth condition on the stochastic differential equations considered. Moreover, results on domain invariance has been proved and the probability of such a preservation has been quantified for the approximate solution given by the NSEM scheme. Due to the fact that the deterministic part is nonstandard, its computed expectation preserves domain invariance unconditionally. We illustrated numerically the main advantages of the NSEM scheme through the geometric Brownian motion compared to the EM scheme and the BIM scheme. \\

The preservation of domain invariance by the NSEM scheme can be used in many applications to select models such as in astronomy (see \cite{bcp}) or biology (see \cite{cps2} and \cite{cps1}). Indeed, as it can be seen in the literature (see \cite{cresson-pierret_nsfdm} and references therein), most of the biological and dynamical population models, simulated with deterministic nonstandard schemes, behave likewise or better than all the classical schemes such as the family of Runge-Kutta methods but without their complex definition and implementation. This is also the case for the balanced implicit methods which need a systematic work to define a proper scheme adapted to the problem considered. \\

Due to the time-dependency and unbounded variations of Brownian motion, further investigations are needed to define more general nonstandard schemes for any stochastic differential equations due to the complexity of the nonstandard rules defined by R.E. Mickens in the stochastic context. This work is in progress and will be the subject of a future paper.

\appendix

\section{Proof of Lemma \ref{lemass1}, \ref{lemass2} and Theorem \ref{convnsem}}

\subsection{Proof of Lemma \ref{lemass1}}
From the continuous-time approximation and the H\"older's inequality \cite[Theorem 189]{hardy} that
\begin{align*}
	|\overline{X}(t)|^{p} &\leq 3^{p-1}\left[|Y_0|^p+\left(\h\right)^p \left|\int_{0}^{t}f(X(s))ds\right|^{p}+\left|\int_{0}^{t}g(X(s))dW(s)\right|^{p}\right] \\
	&\leq 3^{p-1}\left[|Y_0|^{p}+\left(\h\right)^pT^{p-1}\int_{0}^{T}|f(X(s))|^{p}ds+|\int_{0}^{t}g(X(s))dW(s)|^{p}\right] .
\end{align*}
Using the Burkholder-Davis-Gundy \cite[Theorem 7.3]{mao} and H\"older's inequalitiy we obtain for all $0\leq \tau \leq T$,
\begin{align*}
	\E\left[\sup_{0\leq t\leq \tau} |\overline{X}(t)|^{p} \right] \leq 3^{p-1}\bigg[&|Y_0|^{p}+\left(\h\right)^pT^{p-1}\E\left[\int_{0}^{\tau}|f(X(s))|^{p}ds\right] \\
	&+C_p T^{\frac{p}{2}-1}\E\left[\int_{0}^{\tau}|g(X(s))|^{p}ds\right]\bigg]
\end{align*}
where $C_{p}$ is a constant whose expression is given explicitely in function of $p$ in \cite[Theorem 7.3]{mao}. Using the linear growth condition we obtain
\begin{equation*}
	\E\left[\sup_{0\leq t\leq \tau} |\overline{X}(t)|^{p} \right] \leq 3^{p-1}\bigg[|Y_0|^{p}+2^{p-1}K^p\bigg(\left(\h\right)^pT^{p-1}+C_p T^{\frac{p}{2}-1}\bigg)\int_{0}^{\tau}\E\left[\left(1+|X(s)|^p\right)\right]ds\bigg].
\end{equation*}
Then we have
\begin{equation*}
	\E\left[1+\sup_{0\leq t\leq \tau} |\overline{X}(t)|^{p} \right] \leq \tilde{C}_1+\tilde{C}_2 \int_{0}^{\tau}\E\left[1+ \sup_{0\leq r\leq s} |\overline{X}(r)|^{p}\right]ds
\end{equation*}
where $\tilde{C}_1=1+3^{p-1}|Y_0|^{p}$ and $\tilde{C}_2=6^{p-1}K^p\bigg(\left(1+\frac{c(h)}{h}\right)^pT^{p-1}+C_p T^{\frac{p}{2}-1}\bigg)$.\\
Using the Gronwall inequality \cite[Theorem 8.1]{mao} we obtain
\begin{equation*}
	\E\left[\sup_{0\leq t\leq T} |\overline{X}(t)|^{p} \right]\leq \tilde{C}_{1}e^{T \tilde{C}_2}-1 \ .
\end{equation*}
In the same way, we can show that
\begin{equation*}
	\E\left[\sup_{0\leq t\leq T} |\overline{Y}(t)|^{p} \right]\leq \tilde{C}_{1}e^{T \tilde{C}_1}-1.
\end{equation*}
Denoting $C_1=\tilde{C}_{1}e^{T \tilde{C}_2}-1$ we obtain the result.

\subsection{Proof of Lemma \ref{lemass2}}
Let $k_s$ be the integer such that $s\in[t_{k_s}, t_{{k_s}+1}[$. We have
\begin{equation*}
	\left|X(s)-\overline{X}(s)\right|^2 \leq 2\left[\left(\h\right)(s-t_{k_s})^2|f(X_{k_s})|^{2}+|g(X_{k_s})|^{2}|W(s)-W(t_{k_s})|^{2}\right].
\end{equation*}
Now, from the local Lipschitz condition, for $|z|\leq R$ we have
\begin{equation*}
	|f(z)|^{2}\leq 2(|f(z)-f(0)|^{2}+|f(0)|^{2})\leq 2(L_R^2|z|^{2}+|f(0)|^{2}),
\end{equation*}
and, similarly,
\begin{equation*}
	|g(z)|^{2}\leq 2(L_R^2|z|^{2}+|g(0)|^{2}).
\end{equation*}
Then 
\begin{equation*}
	\left|X(s)-\overline{X}(s)\right|^2 \leq 4\left(L_R^2 |X_{k_s}|^2 + |f(0)|^{2}\vee |g(0)|^{2}\right)\left(\varphi(h)^2 +|W(s)-W(t_{k_s})|^{2}\right).
\end{equation*}
Taking the expectation and using Lemma \ref{lemass1} leads to
\begin{equation*}
	\E\left[\left|X(s)-\overline{X}(s)\right|^2\right] \leq 4 \tilde{C}_2\left(\varphi(h)^2 +nh\right)
\end{equation*}
where $\tilde{C}_2=L_R^2 C_1+|f(0)|^{2}\vee |g(0)|^{2}$. Denoting $C_2=4\tilde{C}_2$, we obtain the result.

\subsection{Proof of Theorem \ref{convnsem}}
We define $e(t)=\overline{X}(t)-Y(t)$, $\tau_{R}=\inf\{t\geq 0,|\overline{X}(t)|\geq R\}$, $\rho_{R}=\inf\{t\geq 0,|y(t)|\geq R\}$ and $\theta_{R}=\tau_{R}\wedge\rho_{R}$. \\

Using Young's inequality \cite[Theorem 61]{hardy}, we have for any $\delta>0$ and $p>2$,
\begin{align*}
	\E\left[\sup_{0\leq t\leq T}|e(t)|^{2}\right]&=\E\left[\sup_{0\leq t\leq T}|e(t)|^{2}\mathds{1}_{\{\tau_{R}>T\ and \ \rho_{R}>T\}}\right]+\E\left[\sup_{0\leq t\leq T}|e(t)|^{2}\mathds{1}_{\{\tau_{R}\leq T \ or\ \rho_{R}\leq T\}}\right]\\
	&\leq \E\left[\sup_{0\leq t\leq T}|e(t\wedge\theta_{R})|^{2}\mathds{1}_{\{\theta_{R}>T\}}\right]+\frac{2\delta}{p}\E\left[\sup_{0\leq t\leq T}|e(t)|^{p}\right] + \\ &\quad +\frac{1-\frac{2}{p}}{\delta^{2/(p-2)}}\E\left[\mathds{1}_{\{\tau_{R}\leq T \ or\ \rho_{R}\leq T\}}\right].
\end{align*}
With Lemma \ref{lemass1} we obtain:
\begin{equation*}
	\E\left[\mathds{1}_{\{\tau_{R}\leq T\}}\right]\leq\E\left[\mathds{1}_{\{\tau\leq T\}}\frac{|\overline{X}(\tau_{R})|^{p}}{R^{p}}\right]\leq\frac{1}{R^{p}}\E\left[\sup_{0\leq t\leq T}|\overline{X}(t)|^{p}\right],
\end{equation*}
and in the same way
\begin{equation*}
	\E\left[\mathds{1}_{\{\rho_{R}\leq T\}}\right]\leq\E\left[\mathds{1}_{\{\rho\leq T\}}\frac{|\overline{X}(\rho_{R})|^{p}}{R^{p}}\right]\leq\frac{1}{R^{p}}\E\left[\sup_{0\leq t\leq T}|\overline{X}(t)|^{p}\right]\leq\frac{C_1}{R^{p}}.
\end{equation*}
By definition,
\begin{equation*}
	\E\left[\mathds{1}_{\{\tau_{R}\leq T \ or\ \rho_{R}\leq T\}}\right]=\mathbb{P}(\tau_{R}\leq T\ or \ \rho_{R}\leq T).
\end{equation*}
As 
\begin{equation*}
	\mathbb{P}(\tau_{R}\leq T\ or \ \rho_{R}\leq T)\leq \mathbb{P}(\tau_{R}\leq T)+\mathbb{P}(\rho_{R}\leq T),
\end{equation*}
we obtain
\begin{equation*}
	\E\left[\mathds{1}_{\{\tau_{R}\leq T \ or\ \rho_{R}\leq T\}}\right] \leq\frac{2C_1}{R^{p}}.
\end{equation*}
Also, we have
\begin{equation*}
	\E\left[\sup_{0\leq t\leq T}|e(t)|^{p}\right]\leq 2^{p-1}\E\left[\sup_{0\leq t\leq T}(|\overline{X}(t)|^{p}+|y(t)|^{p})\right]\leq 2^{p}C_1
\end{equation*}
and then we obtain
\begin{equation*}
	\E\left[\sup_{0\leq t\leq T}|e(t)|^{2}\right]\leq \E\left[\sup_{0\leq t\leq T}|\overline{X}(t\wedge\theta_{R})-y(t\wedge\theta_{R})|^{2}\right] +\frac{2^{p+1}\delta C_1}{p}+\frac{2(p-2)C_1}{p\delta^{2/(p-2)}R^{p}}.
\end{equation*}
Now, we will bound the first term of the previous inequality. Using
\begin{equation*}
	Y(t\wedge\theta_{R})=Y_{0}+\int_{0}^{t\wedge\theta_{R}} f(Y(s))ds+\int_{0}^{t\wedge\theta_{R}}g(Y(s))dW(s),
\end{equation*}
and Cauchy-Schwarz's inequality \cite[Theorem 181]{hardy}, we obtain
\begin{align*}
	|\overline{X}(t\wedge\theta_{R})-Y(t\wedge\theta_{R})|^{2}\leq & 3\bigg[T\int_{0}^{t\wedge\theta_{R}}\left|f(X(s))-f(Y(s))\right|^{2}ds \\
	&+\left|\int_{0}^{t\wedge\theta_{R}}g(X(s))-g(Y(s))dW(s)\right|^{2} \\
	&+\left(\frac{c(h)}{h}\right)^2T\int_{0}^{t\wedge\theta_{R}}\left|f(X(s))\right|^{2}ds\bigg].
\end{align*}
Then for all $0\leq \tau \leq T$, we have
\begin{align*}
	\E\left[\sup_{0\leq t\leq \tau}{|\overline{X}(t\wedge\theta_{R})-Y(t\wedge\theta_{R})|^{2}}\right]\leq & 3\Bigg[T\E\left[\int_{0}^{\tau\wedge\theta_{R}}\left|f(X(s))-f(Y(s))\right|^{2}ds\right] \\
	&+\E\left[\sup_{0\leq t\leq\tau}{\left|\int_{0}^{t\wedge\theta_{R}}g(X(s))-g(Y(s))dW(s)\right|^{2}}\right] \\
	&+\tilde{C}_3\left(\frac{c(h)}{h}\right)^2\Bigg].
\end{align*}
where $\tilde{C}_3=2T^2\tilde{C}_2$. Using consecutively Doob's martingal inequality \cite[Theorem 3.8]{mao} and It\^o's isometry \cite[(iii) Theorem 5.8]{mao} we obtain
\begin{align*}
	\E\left[\sup_{0\leq t\leq \tau}{|\overline{X}(t\wedge\theta_{R})-Y(t\wedge\theta_{R})|^{2}}\right]\leq & 3\Bigg[T\E\left[\int_{0}^{\tau\wedge\theta_{R}}\left|f(X(s))-f(Y(s))\right|^{2}ds\right] \\
	&+4\E\left[\int_{0}^{\tau\wedge\theta_{R}}\left|g(X(s))-g(Y(s))\right|^{2}ds\right]\\
	&+\tilde{C}_3\left(\frac{c(h)}{h}\right)^2\Bigg].
\end{align*}
From the local Lipschitz condition we obtain
\begin{align*}
	&\E\left[\sup_{0\leq t\leq \tau}{|\overline{X}(t\wedge\theta_{R})-Y(t\wedge\theta_{R})|^{2}}\right] \\
	&\leq 3L_R^2\left(T+4\right)\E\left[\int_{0}^{\tau\wedge\theta_{R}}\left|X(s)-Y(s)\right|^{2}ds\right]+3\tilde{C}_3\left(\frac{c(h)}{h}\right)^2\\
	&\leq 6L_R^2\left(T+4\right)\E\left[\int_{0}^{\tau\wedge\theta_{R}}\left|X(s)-\overline{X}(s)\right|^{2}+\left|\overline{X}(s)-Y(s)\right|^{2}ds\right]+3\tilde{C}_3\left(\frac{c(h)}{h}\right)^2\\
	&\leq 6L_R^2\left(T+4\right)\Bigg(\E\left[\int_{0}^{\tau\wedge\theta_{R}}\left|X(s)-\overline{X}(s)\right|^{2}ds\right]+\E\left[\int_{0}^{\tau}\sup_{0\leq r\leq s}\left|\overline{X}(r\wedge\theta_{R})-Y(r\wedge\theta_{R})\right|^{2}ds\right]\Bigg)\\
	&\quad +3\tilde{C}_3\left(\frac{c(h)}{h}\right)^2\\
\end{align*}
Using Lemma \ref{lemass2} we obtain
\begin{equation*}
	\E\left[\sup_{0\leq t\leq \tau}{|\overline{X}(t\wedge\theta_{R})-Y(t\wedge\theta_{R})|^{2}}\right] \leq \tilde{C}_5 + \tilde{C}_4\E\left[\int_{0}^{\tau}\sup_{0\leq r\leq s}\left|\overline{X}(r\wedge\theta_{R})-Y(r\wedge\theta_{R})\right|^{2}ds\right]\\
\end{equation*}
where
\begin{align*}
	\tilde{C}_4&=6L_R^2\left(T+4\right), \\
	\tilde{C}_5&=6TL_R^2C_2\left(T+4\right)\left(\varphi(h)^2 +nh\right)+18L_R^2\left(T+4\right)\tilde{C}_2\left(\frac{c(h)}{h}\right)^2.
\end{align*}
Using Gronwall's inequality we obtain
\begin{equation*}
	\E\left[\sup_{0\leq t\leq T}{|\overline{X}(t\wedge\theta_{R})-Y(t\wedge\theta_{R})|^{2}}\right] \leq C_3(h,R).
\end{equation*}
with $C_3(h,R)=\tilde{C}_5 e^{T\tilde{C}_4}$. We finally obtain 
\begin{equation*}
	\E\left[\sup_{0\leq t\leq T}|e(t)|^{2}\right]\leq C_3 +\frac{2^{p+1}\delta C_1}{p}+\frac{2(p-2)C_1}{p\delta^{2/(p-2)}R^{p}}.
\end{equation*}
Let $\epsilon>0$ we can choose $\delta$ such as
\begin{equation*}
	\frac{2^{p+1}\delta C_1}{p}<\frac{\epsilon}{3},
\end{equation*}
then choose $R$ such that 
\begin{equation*}
	\frac{2(p-2)C_1}{p\delta^{2/(p-2)}R^{p}}<\frac{\epsilon}{3},
\end{equation*}
and then choose $h$ sufficiently small such that
\begin{equation*}
	C_3(h,R) <\frac{\epsilon}{3},
\end{equation*}
and then we obtain 
\begin{equation*}
	\E\left[\sup_{0\leq t\leq T}|e(t)|^{2}\right]\leq \epsilon,
\end{equation*}
Using Lyapunov's inequality (\cite[4.12 p.17]{kloeden1}) we obtain
\begin{equation*}
	\E\left[\sup_{0\leq t\leq T}|\overline{X}(t)-Y(t)|\right]\leq\sqrt{\E\left[\sup_{0\leq t\leq T}|\overline{X}(t)-Y(t)|^{2}\right]}\leq \sqrt{\epsilon}
\end{equation*}
that is to say
\begin{equation*}
	\lim_{h\rightarrow 0}\E\left[\sup_{0\leq t\leq T}|\overline{X}(t)-Y(t)|\right]=0.
\end{equation*}
This conclude the proof.

\section{Proof of Proposition \ref{dominv1}, \ref{propcondinv}, \ref{dominv2} and \ref{dominv3}}

\subsection{Proof of Proposition \ref{dominv1}}
\label{dem_dominv1}
Let $I\subset\{1,\dots,n\}$ and let $i\in I$. By hypothesis we have $f_i(a)\ge 0$ and we suppose for $k\ge0$, $(X_k)_i \geq a_i$. Taylor's expansion up to the first order with integral remainder on $f$ gives 
\begin{equation*}
	f_i(X_k) = f_i(a) +\left[\int_0^1 f'(a+s(X_k-a))ds\right]_{ij} (X_k-a)_j .
\end{equation*}
for $1\leq j \leq n$ and where we adopted the Einstein summation convention. Inserting the Taylor's expansion of $f$ in the NSEM gives 
\begin{equation*}
	(X_{k+1}-a)_i = (X_k-a)_j\left(\delta_{ij}+\varphi(h)\left[\int_0^1 f'(a+s(X_k-a))ds\right]_{ij}\right) +\varphi(h)f_i(a)
\end{equation*}
where $\delta_{ij}$ are the components of the identity matrix of $\R^n$. From the hypothesis $f_i(a)\ge 0$ and the definition of $\varphi$, we obtain
\begin{equation*}
	(X_{k+1}-a)_i \ge (X_k-a)_i \phi(Dh).
\end{equation*}
As $(X_k)_i\ge a_i$ and $\phi(Dh)\in]0,1[$ then, $X_{k+1}$ remains in $K^+$ for all $h>0$. This concludes the proof.

\subsection{Proof of Proposition \ref{propcondinv}}
\label{dem_dominv1_nsem}
Let $I\subset\{1,\dots,n\}$ and let $i\in I$. By hypothesis we have $f_i(a)\ge 0$ and $g_{ij}(a)=0$ for $1\leq j \leq d$. We assume for $k\ge0$ that $(X_k)_i \geq a_i$. Taylor's expansion up to the first order with integral remainder on $f$ and $g$ gives 
\begin{align*}
	f_i(X_k) &= f_i(a) + \left[\int_0^1 f'(a+s(X_k-a))ds\right]_{il} (X_k-a)_l, \\
	g_{ij}(X_k) &= g_{ij}(a) + \left[\int_0^1 g'(a+s(X_k-a))ds\right]_{ijl} (X_k-a)_l,
\end{align*}
where $1\leq l \leq n$. Inserting the Taylor's expansion of $f$ and $g$ in the NSEM gives 
\begin{align*}
	(X_{k+1}-a)_i &=\varphi(h)f_i(a) + g_{ij}(a)(\Delta W_k)_j\\
	&+ (X_k-a)_l\bigg(\delta_{il} +\varphi(h) \left[\int_0^1 f'(a+s(X_k-a))ds\right]_{il} \\
	& + (\Delta W_k )_j \left[\int_0^1 g'(a+s(X_k-a))ds\right]_{ijl}\bigg)
\end{align*}
By hypothesis $f_i(a)\ge 0$ and $g_{ij}(a)=0$. Then, using the definition of $\varphi$ we obtain
\begin{equation*}
	(X_{k+1}-a)_i \ge (X_k-a)_i \left(\phi(Dh)-S\left(d\sup_{1\leq p \leq d}|\Delta W_k|_p\right)\right).
\end{equation*}
As $(X_k)_i\ge a_i$ and $\phi(D h)\in]0,1[$ then, $X_{k+1}$ remains in $K^+$ if and only if $h$ is such that
\begin{equation}
	\sup_{1\leq p \leq d}|\Delta W_k|_p\leq \frac{\phi(Dh)}{Sd}.
\end{equation}
This conclude the proof.

\subsection{Proof of Proposition \ref{dominv2}}
\label{dem_dominv2}
Let $k\geq0$ and $1\leq p \leq d$. The probability function is determined by a normal law of zero mean and variance $h$. Using the parity of the normal law, we obtain
\begin{equation*}
	\mathbb{P}\left(|\Delta W_k|_p\leq \frac{\phi(Dh)}{Sd} \right) = 2\mathbb{P}\left((\Delta W_k)_p\leq \frac{\phi(Dh)}{Sd} \right)-1.
\end{equation*}
By definition, we have
\begin{equation*}
	\mathbb{P}\left((\Delta W_k)_p\leq \frac{\phi(Dh)}{Sd} \right) = \int_{-\infty}^{ \frac{\phi(Dh)}{Sd}} \frac{e^{-\frac{x^2}{2h}}}{\sqrt{2\pi h}} dx = \frac{1}{2}\left(1+\mathsf{erf}\left( \frac{\phi(Dh)}{Sd} \frac{1}{\sqrt{2h}}\right)\right)
\end{equation*}
where $\mathsf{erf}$ is the classical error function defined by
\begin{equation*}
	\mathsf{erf}(x)=\frac2{\sqrt\pi} \int_0^x e^{-t^2} dt.
\end{equation*}
Let $0<\epsilon<\frac{1}{2}$. We want to have the probability
\begin{equation*}
	\mathbb{P}\left(\sup_{1\leq p \leq d}|\Delta W_k|_p\leq \frac{\phi(Dh)}{Sd} \right)> 1-\epsilon
\end{equation*}
close to one as long as we chose $\epsilon$ sufficiently small, i.e. by definition we want to have
\begin{equation*}
	2\mathbb{P}\left(\sup_{1\leq p \leq d}(\Delta W_k)_p\leq \frac{\phi(Dh)}{Sd} \right) -1 > 1-\epsilon .
\end{equation*}
Let 
\begin{equation*}
	\frac{1}{2}\left(1+\mathsf{erf}\left( \frac{\phi(Dh)}{Sd} \frac{1}{\sqrt{2h}}\right)\right)=1-\frac{\epsilon}{2}.
\end{equation*}
Using the function $\mathsf{erf}^{-1}$, we obtain
\begin{equation}
	\label{eqh}
	\frac{\phi(Dh)}{Sd\sqrt{2h}}=\mathsf{erf}^{-1}(1-\epsilon).
\end{equation}
Considering $h_0(\epsilon)$ the solution of $\eqref{eqh}$ then, for all $h<h_0(\epsilon)$, we finally obtain
\begin{equation*}
	\mathbb{P}\left(\sup_{1\leq p \leq d}|\Delta W_k|_p\leq \frac{\phi(Dh)}{Sd} \right) \geq \mathbb{P}\left(|\Delta W_k|_p\leq \frac{\phi(Dh)}{Sd} \right) > 1- \epsilon.
\end{equation*}
This concludes the proof.

\subsection{Proof of Proposition \ref{dominv3}}
\label{dem_dominv3}
Let $I\subset\{1,\dots,n\}$ and let $i\in I$. By hypothesis we have $f_i(a)\ge 0$ and we assume for $k\ge0$ that $(X_k)_i \geq a_i$. Taylor's expansion up to the first order with integral remainder on $f$ gives 
\begin{align*}
	f_i(X_k) &= f_i(a) + \left[\int_0^1 f'(a+s(X_k-a))ds\right]_{il} (X_k-a)_l,
\end{align*}
where $1\leq l \leq n$. Inserting the Taylor's expansion of $f$ in the NSEM gives 
\begin{align*}
	(X_{k+1}-a)_i &=\varphi(h)f_i(a) + g_{ij}(X_k)(\Delta W_k)_j\\
	&+ (X_k-a)_l\bigg(\delta_{il} +\varphi(h) \left[\int_0^1 f'(a+s(X_k-a))ds\right]_{il} \bigg)
\end{align*}
We have $\di -D(X_k-a)_i\le(X_k-a)_l \left[\int_0^1 f'(a+s(X_k-a))ds\right]_{il}\leq D(X_k-a)_i$ and $\di\Eb{g(X_k) \cdot \Delta W_k}=0$ by using a classical property of It\^o's integral (see \cite[Theorem 3.2.1 (iii) p. 30]{oksendal}). Then, we obtain
\begin{equation*}
	\Eb{(X_{k+1}-a)_i} = \Eb{(X_k-a)_i}\left(1-\varphi(h)D\right) +\varphi(h)\Eb{f_i(a)}
\end{equation*}
From hypothesis $\Eb{f_i(a)}\ge 0$, $a$ is a constant and the definition of $\varphi$, we obtain
\begin{equation*}
	(\Eb{X_{k+1}}-a)_i \ge (\Eb{X_{k}}-a)_i \phi(D h).
\end{equation*}
As $\Eb{X_k}\ge a$ and $\phi(D h)\in]0,1[$ then, $\Eb{X_{k+1}}$ remains in $K^+$ for all $h>0$. This concludes the proof.

\section{Proof of Lemma \ref{lemmeh0em} and \ref{lemmeh0ns}}

\subsection{Proof of Lemma \ref{lemmeh0em}}
The Euler-Maruyama scheme is the standard version of our scheme which correspond to the case where we have $h<\frac{1}{\lambda}$ and $\phi(\lambda h)=1-\lambda h$. The resolution of \eqref{eqh} gives two solutions
\begin{equation*}
	h_\pm(\epsilon)=\frac{1}{\lambda }  + \frac{\alpha(\epsilon) ^2  \sigma ^2}{\lambda ^2} \pm \frac{ \alpha \sigma  \sqrt{\alpha(\epsilon) ^2  \sigma ^2+2 \lambda }}{\lambda ^2}.
\end{equation*}
The minimal step $h_0(\epsilon)$ is chosen as the minimum between the two roots and $\frac{1}{\lambda}$. As $0<h_-(\epsilon)<h_+(\epsilon)$ and the expression of these roots we obtain $h_0(\epsilon)=h_-(\epsilon)$.

\subsection{Proof of Lemma \ref{lemmeh0ns}}
The resolution of \eqref{eqh} with $\phi(\lambda h)=\exp(-\lambda h)$ leads to solve 
\begin{equation*}
	2\lambda h \exp(2\lambda h)=\frac{\lambda}{\sigma^2 \alpha(\epsilon)^2}.
\end{equation*}
Using the definition of the product logarithm function we obtain the result.


\begin{thebibliography}{17}
	\providecommand{\natexlab}[1]{#1}
	\providecommand{\url}[1]{\texttt{#1}}
	\expandafter\ifx\csname urlstyle\endcsname\relax
	\providecommand{\doi}[1]{doi: #1}\else
	\providecommand{\doi}{doi: \begingroup \urlstyle{rm}\Url}\fi
	
	\bibitem[{Behar} et~al.(2014){Behar}, {Cresson}, and {Pierret}]{bcp}
	E.~{Behar}, J.~{Cresson}, and F.~{Pierret}.
	\newblock {Dynamics of a rotating ellipsoid with a stochastic flattening}.
	\newblock \emph{arXiv}, 1410.0667, 2014.
	
	\bibitem[{Cresson} and {Pierret}(2014)]{cresson-pierret_nsfdm}
	J.~{Cresson} and F.~{Pierret}.
	\newblock {Nonstandard finite difference scheme preserving dynamical
		properties}.
	\newblock \emph{arXiv}, 1410.6661, 2014.
	
	\bibitem[Cresson et~al.(2012)Cresson, Puig, and Sonner]{cps2}
	J.~Cresson, B.~Puig, and S.~Sonner.
	\newblock {S}tochastic models in biology and the invariance problem.
	\newblock \emph{preprint}, 2012.
	
	\bibitem[Cresson et~al.(2013)Cresson, Puig, and Sonner]{cps1}
	J.~Cresson, B.~Puig, and S.~Sonner.
	\newblock {Validating stochastic models: invariance criteria for systems of
		stochastic differential equations and the selection of a stochastic
		Hodgkin-Huxley type model}.
	\newblock \emph{International Journal of Biomathematics and Biostatistics},
	2:\penalty0 111--122, 2013.
	
	\bibitem[Hardy et~al.(1952)Hardy, Littlewood, and P\'{o}lya]{hardy}
	G.H. Hardy, J.E. Littlewood, and G.~P\'{o}lya.
	\newblock \emph{{Inequalities}}.
	\newblock Cambridge university press, 1952.
	
	\bibitem[Higham(2001)]{higham}
	D.J. Higham.
	\newblock {A}n algorithmic introduction to numerical simulation of stochastic
	differential equations.
	\newblock \emph{SIAM review}, 43\penalty0 (3):\penalty0 525--546, 2001.
	
	\bibitem[Hutzenthaler and Jentzen(2015)]{hut_jen}
	M.~Hutzenthaler and A.~Jentzen.
	\newblock \emph{{N}umerical approximations of stochastic differential equations
		with non-globally {L}ipschitz continuous coefficients}, volume 236.
	\newblock American Mathematical Society, 2015.
	
	\bibitem[Kloeden(1994)]{kloeden2}
	P.E. Kloeden.
	\newblock \emph{{Numerical solution of SDE through computer experiments}},
	volume~1.
	\newblock Springer, 1994.
	
	\bibitem[Kloeden and Platen(1992)]{kloeden1}
	P.E. Kloeden and E.~Platen.
	\newblock \emph{{Numerical solution of stochastic differential equations}},
	volume~23.
	\newblock Springer, 1992.
	
	\bibitem[Mao(2007)]{mao}
	X.~Mao.
	\newblock \emph{{Stochastic differential equations and applications}}.
	\newblock Elsevier, 2007.
	
	\bibitem[Mickens(1994)]{mickens1994}
	R.E. Mickens.
	\newblock \emph{{Nonstandard finite difference models of differential
			equations}}.
	\newblock World Scientific, 1994.
	
	\bibitem[Mickens(2005{\natexlab{a}})]{mickens2005}
	R.E. Mickens.
	\newblock {Dynamic consistency: a fundamental principle for constructing
		nonstandard finite difference schemes for differential equations}.
	\newblock \emph{Journal of Difference Equations and Applications}, 11\penalty0
	(7):\penalty0 645--653, 2005{\natexlab{a}}.
	
	\bibitem[Mickens(2005{\natexlab{b}})]{mickens2005adv}
	R.E. Mickens.
	\newblock \emph{{Advances in the Applications of Nonstandard Finite Diffference
			Schemes}}.
	\newblock World Scientific, 2005{\natexlab{b}}.
	
	\bibitem[Milian(1995)]{milian}
	A.~Milian.
	\newblock {Stochastic viability and a comparison theorem}.
	\newblock In \emph{Colloquium Mathematicum}, volume~68, pages 297--316, 1995.
	
	\bibitem[Milstein et~al.(1998)Milstein, Platen, and Schurz]{mil_pla_sch}
	G.N. Milstein, E.~Platen, and H.~Schurz.
	\newblock {Balanced implicit methods for stiff stochastic systems}.
	\newblock \emph{SIAM Journal on Numerical Analysis}, 35\penalty0 (3):\penalty0
	1010--1019, 1998.
	
	\bibitem[{\O}ksendal(2003)]{oksendal}
	B.~{\O}ksendal.
	\newblock \emph{{Stochastic differential equations}}.
	\newblock Springer, 2003.
	
	\bibitem[Schurz(1996)]{schurz}
	H.~Schurz.
	\newblock {Numerical regularization for SDEs: Construction of nonnegative
		solutions}.
	\newblock \emph{Dynamic Systems and Applications}, 5:\penalty0 323--352, 1996.
	
\end{thebibliography}
\end{document}